\documentclass[onecolumn,final,sort&compress]{elsart3p}

\usepackage{comment}
\usepackage{amsmath}       %   AMS equation environments
\usepackage{amssymb}       %   AMS symbol fonts (e.g., \boldsymbol{}
\usepackage{accents}

\usepackage{chapterbib}    %   Bibtex
\usepackage{color,xcolor}
\usepackage{comment}
\usepackage{bm,bbm}
\usepackage{overpic}

\usepackage{times}

\usepackage{soul}

\usepackage{epsfig}
\usepackage{graphicx,graphics,rotating}
\usepackage{subfigure}
\usepackage{multicol}
\usepackage{multirow}

\theoremstyle{plain}

\newtheorem{remark}{Remark}[section]

\definecolor{MyDarkGreen}{rgb}{0,0.45,0}%%{0.80,0.20,0.20} 

\newcommand{\RED}  [1]{{\color{red}  #1}}

\newcommand{\EOD}{\end{document}}

\def\trait #1 #2 #3 {\vrule width #1pt height #2pt depth #3pt}
\def\fin{\hfill
        \trait .3 5 0
        \trait 5 .3 0
        \kern-5pt
        \trait 5 5 -4.7
        \trait 0.3 5 0
\medskip}

\newcounter{numbs}
\newcounter{numbi}
\newcounter{numbii}

\renewcommand{\cv}{\mathbf{c}}

\newcommand{\vv}{\mathbf{v}}

\newcommand{\xv}{\mathbf{x}}

\newcommand{\ds}{d}

\newcommand{\fs}{f}
\newcommand{\gs}{g}
\newcommand{\hs}{h}

\newcommand{\ts}{t}
\newcommand{\us}{u}
\newcommand{\vs}{v}
\newcommand{\ws}{w}
\newcommand{\xs}{x}
\newcommand{\ys}{y}

\newcommand{\Bs}{B}
\newcommand{\Cs}{C}
\newcommand{\Ds}{D}
\newcommand{\Es}{E}

\newcommand{\Hs}{H}

\newcommand{\Ks}{K}
\newcommand{\Ls}{L}

\newcommand{\Ts}{T}

\newcommand{\Ys}{Y}

\newcommand{\calA}{\mathcal{A}}
\newcommand{\calB}{\mathcal{B}}

\newcommand{\calI}{\mathcal{I}}

\newcommand{\calM}{\mathcal{M}}

\newcommand{\REAL}{\mathbbm{R}}
\newcommand{\INTG}{\mathbbm{N}}

\newcommand{\LTWO}{L^2}

\newcommand{\nlen}{\hspace{-0.2mm}}

\newcommand{\norm} [2]{\big|\nlen\big|#1\big|\nlen\big|_{#2}}

\newcommand{\ABS}  [1]{\left|#1\right|}
\newcommand{\abs}  [1]{\big|#1\big|}

\newcommand{\dx}{\,dx}
\newcommand{\dt}{\,dt}
\newcommand{\ddv}{\,dv}

\newcommand{\Omv}{\Omega_{\vs}}
\newcommand{\Omx}{\Omega_{\xs}}

\newcommand{\pp}   {\prime\prime}

\newcommand{\Lst}  {\widetilde{\Ls}}
\newcommand{\hsp}  {\hs^{\prime}}
\newcommand{\hspp} {\hs^{\prime\prime}}
\newcommand{\hsppp}{\hs^{\prime\prime\prime}}

\newcommand{\Hsp}  {\Hs^{\prime}}
\newcommand{\Hspp} {\Hs^{\prime\prime}}

\newcommand{\Dt}{\Delta t}

\newcommand{\expAW}{e^{-\vs^2}}
\newcommand{\expSW}{e^{-\vs^2\slash{2}}}

\newcommand{\calAb}{\overline{\calA}}
\newcommand{\calBb}{\overline{\calB}}

\newcommand{\Csb} {\Cs^\star}

\newcommand{\bigdot}[1]{%
  \accentset{\mbox{\large\bfseries .}}{#1}}

\usepackage{showlabels}
\usepackage[color,notref,notcite]{showkeys}

\begin{document}

\begin{frontmatter}

  \title{Stability and conservation properties of Hermite-based
    approximations of the Vlasov-Poisson system}

  \author[Unimore,CNR]{D. Funaro}
  \author[LANLT5]     {and G. Manzini}
  %\author[LANLT5]     {, and G.~L. Delzanno}

  \address[Unimore]{
    Dipartimento di Scienze Chimiche e Geologiche,
    Universit\`a degli Studi di Modena \\ e Reggio Emilia,
    Italy;
    \emph{e-mail: daniele.funaro@unimore.it}
  }
  \address[CNR]{
    Istituto di Matematica Applicata e Tecnologie Informatiche,
    Consiglio Nazionale \\delle Ricerche,
    via Ferrata 1,
    27100 Pavia,
    %\emph{e-mail: daniele.funaro@unimore.it}
  }
  \address[LANLT5]{
    Group T-5,
    Applied Mathematics and Plasma Physics,
    Theoretical Division,\\
    Los Alamos National Laboratory,
    Los Alamos, NM,
    USA;
    %\emph{e-mail: gmanzini@lanl.gov, delzanno@lanl.gov}
    \emph{e-mail: gmanzini@lanl.gov}
  }
 
  \maketitle
  \begin{abstract}
    Spectral approximation based on Hermite-Fourier expansion of the
    Vlasov-Poisson model for a collisionless plasma in the
    electrostatic limit is provided, by including high-order
    artificial collision operators of Lenard-Bernstein type.
    These differential operators are suitably designed in order to
    preserve the physically-meaningful invariants (number of
    particles, momentum, energy).
    In view of time-discretization, stability results in appropriate
    norms are presented. In this study, necessary conditions link the
    magnitude of the artificial collision term, the number of spectral
    modes of the discretization, as well as the time-step.
    The analysis, carried out in full for the Hermite discretization
    of a simple linear problem in one-dimension, is then partly
    extended to cover the complete nonlinear Vlasov-Poisson model.
  \end{abstract}

\end{frontmatter}
 
\raggedbottom

%% ----------------
%% PAPER
%% ----------------

%% \input{intro_new.tex}
%% \input{hermite_paper.tex}
%% \input{conclusion.tex}

\section{Introduction}
The numerical approximation of physical systems described by kinetic
equations is a formidable challenge~\cite{Pareschi-Dimarco:2014}.
These equations are, indeed, highly dimensional, strongly non-linear,
and describe phenomena that are extremely multi-scale, as the behavior
of the physical system at macroscopic scales is influenced by the
microscopic particle dynamics.
In plasma physics, scale separation occurs at the kinetic level
because of the difference in mass between electrons and
ions~\cite{Hazeltine-Waelbroeck:1998}.
%, and {\color{red}{it}} is exacerbated when one compares 
%{\color{red}{system and kinetic scales - non 
% capisco}}~\cite{Bittencourt:2004}.
%%
Other important applications that may be worth mentioning can be found
in fluid dynamics, particularly, atmospheric and climate
research~\cite{Arakawa-Jung:2011}, and multidimensional radiative transfer
problems~\cite{Luo-Wang-Zhang-Yi-Tan:2018}.
%% and cosmology.\textbf{[Add citation]}
%%
In all these fields, performing macroscale simulations that accurately
include effects from the underlying microscale particle dynamics is
still an open challenge.

\smallskip
In this work, we focus on the numerical approximation of the kinetics
equation describing the behavior of electrically charged particle in a
noncollisional plasmas, also known as \textit{the Vlasov equation}.
Such equation governs the time evolution of the distribution function
of the plasma particles, through the action of an electromagnetic
field generated by the charge and current densities of the same moving
particles.
The resulting coupling through Maxwell's equations (or the Poisson's
equation in the electrostatic limit) is highly nonlinear since the
electromagnetic sources in such equations, i.e., charge and current
densities, depend on the same distribution
functions~\cite{Glassey:1996}.

\smallskip
In his historical and pioneering paper, cf.~\cite{Grad:1949}, Grad
proposed to expand the velocity distribution function of a
noncollisional plasma at equilibrium using Hermite functions.
Hermite functions are Hermite polynomials multiplied by the Gaussian
exponential function, $\ws(\vs)=\exp(-\vs^2)$, where $\vs$ is the
velocity of the plasma particles.
Such a weight $\ws$ is indeed the velocity distribution of a plasma at
equilibrium and is a steady state solution of the Vlasov equation.
Since a plasma at equilibrium is described by the first mode of the
Hermite expansion, we expect that only a few modes may be needed to
describe a plasma in a perturbed state but still close to the
equilibrium.
Moreover, when the solution of the Vlasov equation is expanded on the
Hermite basis functions, the equations for the first three
coefficients correspond to the conservation laws for the number of
particles, momentum and energy, and determine the macroscopic (i.e.,
fluid) behavior of a plasma.
The following terms of the Hermite expansion introduces kinetics
effects in the model in a very straighforward manner, thus providing a
strategy to realize the coupling between micro- and macro-physics.
Thus, the micro/macro coupling is an intrinsic and specific feature of
the Hermite approach, which cannot be replicated if we choose a
different set of basis functions.
For these reasons, Hermite functions are a sort of ``ideal'' basis for
solving numerically Vlasov-based models of noncollisional plasmas.

\smallskip

Since late sixties throughout the last five decades, Grad's idea has
extensively been applied to the development of plasma simulators; see,
for example,
\cite{Armstrong:1970,%
  Gajewski-Zacharias:1977,%
  Klimas:1983,%
  Holloway:1996,%
  Schumer-Holloway:1998,%
  Camporeale-Delzanno-Lapenta-Daughton:2006,%
  Parker-Dellar:2015},
where the Hermite basis for velocity is coupled with the Fourier basis
in space.
A renewed interest has been manifested in very recent years towards
these approximation
methods~\cite{Camporeale-Delzanno-Bergen-Moulton:2015,%
  Camporeale-Delzanno-Bergen-Moulton:2016,%
  Delzanno:2015,Delzanno-Roytershtein:2018}, as the excellent
properties mentioned above make them the natural numerical framework
of high resolution and computationally efficient
solvers~\cite{Vencels-Delzanno-Manzini-Markidis-BoPeng-Roytershteyn:2015,%
  Roytershtein-Boldyrev-Delzanno-Chen-Groselij-Loureiro:2018}.
%%
%% \RED{\textbf{Introduzione modificata da qui.}}\\ \noindent
%%
Moreover, the accuracy of Hermite's approximations can be improved by
order of magnitudes by introducing a \textit{translation factor},
$\us$, and a \textit{scaling factor}, $\alpha$, in the so-called
\textit{generalized weight},
$\ws(\vs)=exp\big(-((\vs-\us)\slash{\alpha})^2\big)$,
cf.~\cite{Tang:1993}.
Empirical evidence that a convenient choice of the scaling factor
$\alpha$ can improve the accuracy in Hermite discretizations of the
Vlasov equation was shown in~\cite{Schumer-Holloway:1998}.
Generalized basis function of Hermite type has been investigated for
solving time-dependent parabolic problems in~\cite{Ma-Sun-Tang:2005}
and, more recently, in~\cite{Fatone-Funaro-Manzini:2018} for the
approximation of the Vlasov phase space.
An adaptive strategy is currently under investigation,
see~\cite{Pagliantini-Delzanno-Manzini-Markidis:2019}, where both
$\us$ and $\alpha$ may change through momentum and energy following
how the plasma evolves in time during a numerical simulation.
Such adaptive strategy is sought to improve the computational
efficiency by using only a few spectral modes where a macroscopic
description of the system is appropriate and adding more modes where
the microscopic physics is
important~\cite{Vencels-Delzanno-Johnson-BoPeng-Laure-Markidis:2015}.
This aspect offers the possibility of selecting the most meaningful
number of spectral modes for a given resolution in phase space.

\indent
The strong point in favour of spectral schemes is that such schemes
can be extremely accurate because of their exceptional convergence
rate, see, for example, the books referenced
in~\cite{Canuto-Hussaini-Quarteroni-Zang:1988,%
  Canuto-Hussaini-Quarteroni-Zang:2006,%
  Bernardi-Maday:1997,%
  Funaro:1992,%
  Funaro:1997,%
  Shen-Tang-Wang:2011}.
Their stability for Vlasov-based systems can be ensured in different
ways.
If we assume that the velocity domain remains bounded during a plasma
simulation, we can use the different basis provided by Legendre
polynomials and stability can be enforced somehow through a penalty
technique acting on boundary terms, see for
example~\cite{Manzini-Delzanno-Vencels-Markidis:2016,Manzini-Funaro-Delzanno:2016}.
Relaxing this assumption yields an unbounded velocity domain and this
approach is no longer feasible.
In the more general case, the Vlasov equation describe a collisionless
transport phenomenon in a six-dimensional phase space, and a
straightforward way to enforce numerical stability to the
discretization of an advection equation is by adding a suitable
artificial dissipation to its, otherwise zero, right-hand side.
However, in the case of the Vlasov equation, using an artificial
dissipation term introduces a major issue because such modification
must not destroy the conservation properties of the original method.
Discrete analogs of the total number of particles (also proportional
to mass and charge of the plasma particles), the total momentum and
the total energy may indeed exist in spectral-based discretizations
using, for the space term, the Fourier expansion
\cite{Holloway:1996,Schumer-Holloway:1998,%
  Camporeale-Delzanno-Bergen-Moulton:2016,%
  Delzanno:2015}, or the discontinuous Galerkin method \cite{%
  Manzini-Delzanno:2017a,Manzini-Delzanno:2017b,%
  Manzini-Koshkarov-Delzanno:2019,
  Koshkarov-Manzini-Delzanno-Pagliantini-Roytershtein:2019:ConsLaws,%%
  Koshkarov-Manzini-Delzanno-Pagliantini-Roytershtein:2019%
}
Conservation properties are fundamental in long-time integration runs
since they provide physically meaningful constraints on the numerical
approximation of the plasma behavior.
Such constraints are strongly related to significant properties like
the well-posedness and robustness of the method, and the reliability
of the numerical simulation.
This fact justifies the great effort that has been devoted in design
spectral methods with such discrete conservation properties.

In the spectral discretizations of the Vlasov equation using Hermite
basis functions, the conservation of number of particles, momentum and
energy is strictly related to the lowest-order modes and can be
destroyed by the numerical dissipation term.
A possible way to maintain a perfect preservation of low modes, is to
design such dissipation terms through Lenard-Bernstein-like operators
(see~\cite{Lenard-Bernstein:1958}) of order $2k$, with integer $k\geq
1$.
In this case, the $1D-1D$ Vlasov-Poisson system of equations takes the
form
\begin{align}
  &\frac{\partial\fs}{\partial t} + \vs\frac{\partial\fs}{\partial x} 
  -\Es\frac{\partial\fs}{\partial\vs} = 
  -(-1)^k \nu\Lst^{(k)}\Ls^{(k)}\fs 
  \qquad\textrm{in~}\Omega\times[0,T],
    \label{eq:Vlasov:1}
  \\[0.5em]
  &\frac{\partial\Es}{\partial x} = 1 - \int_{\Omv}\fs\ddv 
  \qquad\textrm{in~}\Omega\times[0,T],
  \label{eq:Vlasov:2}
\end{align}
where $f$ is the distribution function, $\Es$ the electric field,
%that is consistently generated by the electric charge associated
%with $\fs$; 
$\Lst^{(k)}$ and $\Ls^{(k)}$ are the Lenard-Bernstein-like operators
only acting onto the velocity variable $\vs$. 
The positive parameter $\nu$ is a sort of \emph{artificial viscosity}
used to tune the action of the differential operator
$\Lst^{(k)}\Ls^{(k)}$ on $\fs$.
The combination $\Lst^{(k)}\Ls^{(k)}$ is the Lenard-Bernstein-like
operator of order $2k$ and introduces a sort of \textit{artificial
  collisional term}, i.e., a numerical dissipation, in the equation.
This kind of dissipation terms were proposed in previous works to
control the filamentation process based on an empirical argument, cf.~
\cite{Camporeale-Delzanno-Bergen-Moulton:2015,%
  Camporeale-Delzanno-Bergen-Moulton:2016,%
  Delzanno:2015,Delzanno-Roytershtein:2018,%
  Manzini-Delzanno-Vencels-Markidis:2016,%
  Manzini-Funaro-Delzanno:2016}.

Commonly, there are two different choices of Hermite functions, which
are Hermite polynomials multiplied by a suitable weight function.  
The classical polynomial orthogonality weighted by $\ws(\vs)=\expAW$
leads to the so called \emph{asymmetrically weighted} (AW) case,
whereas the orthogonality of Hermite functions, each one weighted by
$\ws(\vs)=\expSW$, leads to the \emph{symmetrically weighted} (SW)
case.
This terminology will be better clarified in the coming sections.
Accordingly, we have two different definitions of the Lenard-Bernstein
differential operators $\Lst^{(k)}$ and $\Ls^{(k)}$.
%%
%% Once chosen a set of basis functions, we obtain the discrete weak
%% form of the Vlasov equation~\eqref{eq:Vlasov:1} by expanding the
%% distribution function $\fs$ on such Hermite functions and, then,
%% testing the resulting equation on the corresponding set of dual
%% functions.
%%
%Nevertheless, two crucial properties are common to both formulations.
%%
%First, 
%their combination, i.e., $\Lst^{(k)}\Ls^{(k)}$ in the
%right-hand side of~\eqref{eq:Vlasov:1}, is a $2k$-order diffusive
%differential operator.
%%
In both cases, the basis elements are eigenfunctions of the combined
operator, and the corresponding eigenvalues are zero regarding the
first $k-1$ modes.
This actually says that the action of these operators does not modify
such modes, or, in other words, $\Lst^{(k)}\Ls^{(k)}$ induces
dissipation only for the modes starting from $k$. %% $k\geq1$.
%%
%This property implies the conservation of all the moments of $\fs$, up to a certain degree,
%with respect to the velocity variable $\vs$, also in the diffusive regime. 
%%
Despite these common properties, the two discrete formulations
resulting from using AW and SW Hermite functions are substantially
different.
In fact, it turns out that, concerning time-discretization, the SW
formulation can easily be proven to be algebraically stable with or
without the diffusive term
(see~\cite{Holloway:1996,Schumer-Holloway:1998}), while for the AW
formulation the issue is far more delicate.
More precisely, the stability result in the $L^2(\Omega)$ norm that we
are interested to investigate reads as
\begin{align*}
  \frac{d}{\dt}\norm{\fs(\cdot,\cdot,t)}{L^2(\Omega)}^2 \leq 0.
\end{align*}
This inequality trivially implies the boundedness in time of $\fs$.
%$\norm{\fs(\cdot,\cdot,t)}{L^2(\Omega)}\leq\norm{\fsz}{L^2(\Omega)}$
%for every $t\in [0,T]$.
%%
The main criticism to the SW formulation is that, although stable, it
does not effectively preserve the lowest modes during time evolution.
On the contrary, the AW formulation perfectly conserves all the basic
invariants, but its stability needs a deeper analysis. 
What we are able to prove in our work is an $L^2(\Omega)$ stability
result when $\nu$ is sufficiently large thanks to a suitable extension
of the Poincar\`e inequality in weighted norms defined on the real
line.
%% 
%%Despite the efforts devoted to this problem in the past, we were able
%%to succeed
The property of stability then follows by classical estimates for
bilinear forms in Sobolev spaces.  
When instead $\nu$ is small, the result is certainly not true in the
continuous case, but still holds in the framework of numerical
discretizations, by suitably linking $\nu$ to the time discretization
parameter $\Delta t$, the final time $T$, and the maximum integer $N$
used for the Hermite truncation in the variable $v$.
We show how to get these relations for a simple linear
advection-diffusion model problem, and successively we partly extend
our arguments to equation~\eqref{eq:Vlasov:1}.

%Although this property is always satisfied by the exact solution, it
%is often not true for the numerical approximation.
%%
%For example, stability is an intrinsic property of the numerical
%formulation using the SW Hermite basis
%functions~\cite{Holloway:1996,Schumer-Holloway:1998}.
%%
%Here, the combined Lenard-Bernstein diffusive operator
%in~\eqref{eq:Vlasov:1} is normally used to take under control
%nonlinear numerical instabilities such as the filamentation.
%%
%\RED{parti da rivedere ---
%In other situations, it is enforced somehow, for example, %through a
%penalty technique acting on boundary terms when using the basis
%provided by Legendre polynomials as
%in~\cite{Manzini-Delzanno-Vencels-Markidis:2016,Manzini-Funaro-Delzanno:2016},
%or by adding special artificial diffusion terms, as is done when %using
%the AW Hermite basis functions.
%%
%Instead, when using AW Hermite basis functions, the %Lenard-Bernstein
%operator is really needed to control the stability of the %scheme.}

%In the tourough investigation that we carry out in this paper, we
%derive a set of inequalities relating $\nu$, the timestep $\Delta$,
%the total number of Hermite modes $N$, and the final time of a
%numerical simulation $T$.
%%
%Such conditions are sufficient to guarantee the numerical stability of
%the AW method for the (simple) advection equation and for the Vlasov
%equation.
%%
%According to our results, stability in the $\LTWO$-norm can only be
%provided when $\nu\geq1$.
%%
%However, for values of $\nu$ smaller than $1$, a stability result is
%still obtained in a weaker norm.

%% 
A stability result for the Hermite approximation
of $1D-1V$ Vlasov-Poisson model was provided
in~\cite{Gajewski-Zacharias:1977}, where $\LTWO$ boundedness is
proven with respect to the parameter $N$.
However, that paper fails in proving absolute stability with
respect to $t$, since the estimate there provided contains an
exponential growth in time on the right-hand side of the estimate
inequality.
The major result of our work is in achieving a stability estimate
where boundedness in time is guaranteed for all $t$.

\smallskip
The outline of the paper is as follows.
In Section~\ref{sec:preliminary:properties}, we introduce the
Hermite-based discretization framework and discuss some useful
relations.
In Sections~\ref{sec:diffusive:operators:AW}
and~\ref{sec:action:operators:AW}, we introduce the
Lenard-Bernstein-like operators for the spectral method using the
asymmetrically weighted (AW) Hermite functions, and study their
actions on the conservation property of the Vlasov-Poisson system.
In Sections~\ref{sec:diffusive:operators:SW}
and~\ref{sec:action:operators:SW}, we do the same for the spectral
method using the symmetrically weighted (SW) Hermite functions.
In Section~\ref{sec:Hermite:advection} we introduce the SW and AW
Hermite discretization of the advection problem
\begin{align}
  \frac{\partial\fs}{\partial t} - \frac{\partial\fs}{\partial\vs} = -(-1)^k \nu\Lst^{(k)}\Ls^{(k)}\fs,
  \label{eq:model:1}
\end{align}
for the unknown scalar field $\fs(\vs,t)$, with the initial condition
$\fs(\vs,0)=\fs_0(\vs)$, and in Section~\ref{sec:adve:stab}, we study
how the stabilization operator $\Lst^{(k)}\Ls^{(k)}$ impact on its
spectral discretization.
In Section~\ref{sec:time:discretization}, we apply the implicit time
discretization to the system of coefficient resulting from the Hermite
discretization and investigate its stability using a suitable weighted
norm.
In Section~\ref{sec:full-discretization:Vlasov-Poisson}, we extend our
approach to the full spectral discretization of the Vlasov-Poisson
system of equations, and derive sufficient condition to guarantee the
stability of the method.
In Section~\ref{sec:conclusion} we offer our final remarks and
conclusions.

\section{Preliminary properties of the Hermite polynomials}
\label{sec:preliminary:properties}

\smallskip
%$\REAL=\REAL$, $\Omx=[0,2\pi]$ or $\Omx=[-\pi,\pi]$

\noindent 
We start by pointing out some well-known relations concerning Hermite
polynomials, that, as usual, are denoted by $\Hs_{n}(\vs)$ and we
consider as functions of the independent variable $\vs\in\REAL$, the
integer number $n$ being the degree of the polynomial.
First of all, we have the three-point recursion formula that links
$\Hs_{n+1}$ to $\Hs_{n}$ and $\Hs_{n-1}$:
\begin{align}
  &\Hs_0=1,\quad\Hs_1=2\vs, \\[0.5em]
  &\Hs_{n+1} = 2\vs\Hs_{n} - 2n\Hs_{n-1},
  \qquad n\geq 1
\end{align}
and the differential equation for $\Hs_{n}$ 
\begin{align}
  \Hs_{n}^{\pp} - 2\vs\Hs_{n}^{\prime} + 2n\Hs_{n} = 0, %%, \qquad \forall n\in\INTG.
  \label{eq:0075}
\end{align}
which holds for $n\in\INTG$ and where $\prime$ and $\pp$ denote the
first and second derivatives with respect to $v$.
Moreover, the next formulas link Hermite polynomials of different
degrees $n$:
\begin{align}
  \Hs_{n}^{\prime} = 2\vs\Hs_{n} - \Hs_{n+1},% \qquad \forall n\in\INTG,
  \label{eq:0hermite:2}
\end{align}
\begin{align}
  \Hs^\prime_0=0 \quad {\rm and}\quad \Hs_{n}^{\prime} = 2n\Hs_{n-1}, \qquad \forall n\geq 1.
  \label{eq:0hermite:3}
\end{align}
The relation between the Hermite polynomials and their first
derivative in~\eqref{eq:0hermite:3} can recursively be generalized as
follows:
\begin{align}
  \Hs^{(m)}_n =
  \begin{cases}
    0                             & \mbox{$n<m$},\\
    2^m\frac{n!}{(n-m)!}\Hs_{n-m} & \mbox{$n\geq m$}.
  \end{cases}
  \label{eq:0hermite:3a}
\end{align}

Hermite polynomials are orthogonal with respect to the weight function
$\expAW$ and are normalized in such a way that:
\begin{align}
  \int_{\REAL}\Hs^2_{n}\expAW\ddv = \sqrt{\pi}\,2^{n}\,n!.
  \label{eq:Hermite:polynomial:norm}
\end{align}

By examining relation~\eqref{eq:0hermite:3}, it turns out that the
derivatives of the Hermite polynomials are also orthogonal with
respect to the weight $\expAW$.
Using~\eqref{eq:0hermite:3} and~\eqref{eq:Hermite:polynomial:norm} for
$n\geq 1$, we can find that:
\begin{align}
  \int_{\REAL}\big(\Hsp_n\big)^2\expAW\ddv
  &= 4n^2\int_{\REAL} \big(\Hs_{n-1}\big)^2\expAW\ddv
  = 4n^2\,\sqrt{\pi}\,2^{n-1}\,(n-1)! \nonumber \\[0.5em]
  &= 2n\,\sqrt{\pi}\,2^{n}\,n! 
  = 2n\int_{\REAL} \Hs^2_{n}\expAW \ddv .
  \label{eq:005}
\end{align} 
The above relation is trivially satisfied also for $n=0$.
For $n>m$, we recursively find that
\begin{align}
  \int_{\REAL}\big(\Hs^{(m)}_n\big)^2\expAW\ddv
  = 2^{m}\frac{n!}{(n-m)!}\int_{\REAL}\Hs^2_{n}\expAW \ddv.
  \label{eq:005a}
\end{align} 

\smallskip
Consider the generic function $\varphi$ that can be expanded as a
series of Hermite polynomials
$\varphi=\sum_{n=0}^\infty\Cs_{n}\Hs_{n}$ and its first derivative
$\varphi^{\prime}=\sum_{n=1}^\infty\Cs_{n}\Hsp_{n}$.
The Fourier coefficients $\Cs_{n}$ of $\varphi$ are obtained as usual:
\begin{align}
  \Cs_n =\frac{1}{\sqrt{\pi}\,2^{n}\,n!}\int_{\REAL}\varphi \Hs_n \expAW\ddv. %%,\quad n\geq 0.
  \label{eq:0052}
\end{align} 
Of course, $\varphi$ has to be such that all the above integrals are
finite.
% and we suppose that can be expanded in the Hermite basis in an
% appropriate way.
%%
From the orthogonality of Hermite polynomials and their
derivatives, it follows that:
\begin{align*}
  \int_{\REAL}\varphi^2   \expAW\ddv &= \sum_{n=0}^{\infty}\Cs_{n}^2\int_{\REAL}\Hs_{n}^2   \expAW\ddv,\\[0.5em]
  \int_{\REAL}(\varphi')^2\expAW\ddv &= \sum_{n=0}^{\infty}\Cs_{n}^2\int_{\REAL}(\Hs_{n}')^2\expAW\ddv.
\end{align*}
The last summation can also start from $n=1$ since $\Hs_{0}'=0$.

\smallskip
We show a few inequalities that will be used later in this paper.
%%
%% We begin with proving a Poincar\'e-type inequality.
%% \medskip\noindent
%%
By isolating the effect of the first Fourier coefficients, we can
prove Poincar\'e-type inequalities for a linear combination of Hermite
polynomials and their first derivatives with respect to the norm
induced by the weighted $L^2$ inner product where the weight is equal
to $\expAW$.
%%
% Consider the polynomial function
% $\varphi=\sum_{n=0}^\infty\Cs_{n}\Hs_{n}$ and its first derivative
% $\varphi^{\prime}=\sum_{n=1}^\infty\Cs_{n}\Hsp_{n}$.
%%
%%
Indeed, the orthogonality of the first derivatives of the Hermite
polynomials, equation~\eqref{eq:005}, and the fact that $2n\geq2$ for
$n\geq1$, imply that:
\begin{align}
  \int_{\REAL}\big(\varphi^{\prime}\big)^2\expAW\ddv
  &
  =\int_{\REAL}\big( \sum_{n=1}^\infty\Cs_{n}\Hsp_{n} \big)^2\expAW\ddv
  = \sum_{n=1}^\infty\Cs_{n}^2\int_{\REAL}\big(\Hsp_{n}\big)^2\expAW\ddv\nonumber\\[0.5em]
  &
  = \sum_{n=1}^\infty\Cs_{n}^2\,2n\,\int_{\REAL}\Hs^2_{n}\expAW\ddv
  \geq 2\sum_{n=1}^\infty\Cs_{n}^2\int_{\REAL}\Hs_{n}^2\expAW\ddv,
  \label{eq:2501}
\end{align}
where all summations start from $n=1$ since $\Hs_{0}=1$ and
$\Hsp_{0}=0$.
Then, we add and subtract the weighted integral of the zeroth-order mode, i.e,
$\Cs^2_{0}\Hs^2_{0}$, to the last member of inequality~\eqref{eq:2501}
and use the expansion of $\varphi$, so to have
\begin{align}
  \int_{\REAL}\big(\varphi^{\prime}\big)^2\expAW\ddv
  % &\geq 
  % 2\left(
  %   \sum_{n=1}^\infty\Cs_{n}^2\int_{\REAL}\Hs_{n}^2\expAW\ddv
  %   +\Cs_{0}^2\int_{\REAL}\Hs_{0}^2\expAW\ddv-\Cs_{0}^2\int_{\REAL}\Hs_{0}^2\expAW\ddv
  % \right)
  % \nonumber\\[0.5em]
  &\geq %&=
  2\sum_{n=0}^\infty\Cs_{n}^2\int_{\REAL}\Hs_{n}^2\expAW\ddv
  -2\Cs_{0}^2\int_{\REAL}\Hs_{0}^2\expAW\ddv
  \nonumber\\[0.5em]
  &= 2\int_{\REAL}\varphi^2\expAW\ddv - 2 \sqrt{\pi}\Cs_{0}^2.
  \label{eq:2502}
\end{align}
By reversing this inequality we find that
\begin{align}
  \int_{\REAL}\varphi^2\expAW\ddv
  \leq
  \frac{1}{2}\int_{\REAL}\big(\varphi^{\prime}\big)^2\expAW\ddv
  +\sqrt{\pi}\Cs_{0}^2.
  \label{eq:2503}
\end{align}

This inequality can be generalized to derivatives of order $m>1$.
Since $\Hs^{(m)}_{n}=0$ for $n<m$, using
formulas~\eqref{eq:0hermite:3a} and~\eqref{eq:005a}, we find that
\begin{align}
  \int_{\REAL}\big(\varphi^{(m)}\big)^2\expAW\ddv
  &
  = \int_{\REAL}\big( \sum_{n=m}^\infty\Cs_{n}\Hs^{(m)}_{n} \big)^2\expAW\ddv
  = \sum_{n=m}^\infty\Cs_{n}^2\int_{\REAL}\big(\Hs^{(m)}_{n}\big)^2\expAW\ddv\nonumber\\[0.5em]
  &
  = \sum_{n=m}^\infty\Cs_{n}^2\,2^{m}\frac{n!}{(n-m)!}\,\int_{\REAL}\Hs^2_{n}\expAW\ddv
  \geq 2^{m}\,m!\sum_{n=m}^\infty\Cs_{n}^2\int_{\REAL}\Hs_{n}^2\expAW\ddv,
  \label{eq:2504}
\end{align}
as $n!\slash{(n-m)!}\geq m!$ when $n\geq m$.
Now, we add and subtract the weighted integral of the first $m$ modes,
i.e., $\big(\Cs_{\ell}\Hs_{\ell}\big)^2$, $\ell=0,\ldots,m-1$, to the
last member of~\eqref{eq:2504}, and use the normalization of the
Hermite polynomials to find that
\begin{align}
  \int_{\REAL}\big(\varphi^{(m)}\big)^2\expAW\ddv
  % &\geq 2^{m}\left(
  %   \sum_{n=m}^\infty\Cs_{n}^2\int_{\REAL}\Hs_{n}^2\expAW\ddv
  %   +\sum_{\ell=0}^{m-1}\Cs_{\ell}^2\int_{\REAL}\Hs_{\ell}^2\expAW\ddv
  %   -\sum_{\ell=0}^{m-1}\Cs_{\ell}^2\int_{\REAL}\Hs_{\ell}^2\expAW\ddv
  % \right)
  % \nonumber\\[0.5em]
  &\geq 
  2^{m}\,m!\left(
    \sum_{n=0}^\infty\Cs_{n}^2\int_{\REAL}\Hs_{n}^2\expAW\ddv
    -\sum_{\ell=0}^{m-1}\Cs_{\ell}^2\int_{\REAL}\Hs_{\ell}^2\expAW\ddv
  \right)
  \nonumber\\[0.5em]
  &= 2^{m}\,m!\left(
    \int_{\REAL}\varphi^2\expAW\ddv
    -\sqrt{\pi}\sum_{\ell=0}^{m-1}\,2^{\ell}\,\ell!\,\Cs_{\ell}^2
  \right).
  \label{eq:2506}
\end{align}
By reversing this inequality we find that
\begin{align}
  \int_{\REAL}\varphi^2\expAW\ddv
  \leq 
  \frac{1}{2^{m}\,m!}\int_{\REAL}\big(\varphi^{(m)}\big)^2\expAW\ddv
    +\sqrt{\pi}\sum_{\ell=0}^{m-1}\,2^{\ell}\,\ell!\,\Cs_{\ell}^2.
  \label{eq:2507}
\end{align}
The most general Poincar\'e-type inequality is the one involving
derivatives of order $m$ and $p$.
%%
%%Assume, with no loss of generality, that $m>p$. 
Assuming that $m>p$ and noting that
$\Hs^{m}=\Hs^{(p+(m-p))}=\big(\Hs^{(p)}\big)^{(m-p)}$, a
straightforward calculation exploiting the orthogonality of the
derivatives of the Hermite polynomials yields
\begin{align}
  \int_{\REAL}\big(\varphi^{(m)}\big)^2\expAW\ddv 
  &
  = \int_{\REAL}\big( \sum_{n=m}^\infty\Cs_{n}\Hs^{(m)}_{n} \big)^2\expAW\ddv
  = \int_{\REAL}\big( \sum_{n=m}^\infty\Cs_{n}\Hs^{(p+(m-p))}_{n} \big)^2\expAW\ddv
  \nonumber\\[0.5em]
  &= \sum_{n=m}^\infty\Cs_{n}^2\,2^{m-p}\frac{n!}{(n-(m-p))!}\,\int_{\REAL}\big(\Hs^{(p)}_{n}\big)^2\expAW\ddv
  \nonumber\\[0.5em]
  &\geq 2^{m-p}\,(m-p)!\sum_{n=m}^\infty\Cs_{n}^2\int_{\REAL}\big(\Hs_{n}^{(p)}\big)^2\expAW\ddv,
  \label{eq:2506a}
\end{align}
where we also used the fact that $n!\slash{(n-(m-p))!}>(m-p)!$ for $n>1$.
Then, we add and subtract the weighted integrals of 
%% the $(m-p)$ intermediate (derivative) modes, i.e.,
$\Cs_{\ell}^2\big(\Hs_{\ell}^{(p)}\big)^2$
for $\ell=p,\ldots,p+(m-p)-1$, to the last member of~\eqref{eq:2506a}
and we repeat the same argument as above to obtain 
%% find that
% \begin{align}
%   \int_{\REAL}\big(\varphi^{(m)}\big)^2\expAW\ddv
%   &\geq 2^{m-p}
%   \left(
%     \sum_{n=m}^\infty  \Cs_{n}^2  \int_{\REAL}\big(\Hs_{n}^{(p)}   \big)^2\expAW\ddv +
%     \sum_{\ell=p}^{m-1}\Cs_{\ell}^2\int_{\REAL}\big(\Hs_{\ell}^{(p)}\big)^2\expAW\ddv
%   \right.
%   \nonumber\\[0.5em]
%     &
%     \left.
%     \phantom{
%       \geq 2^{m-p}\left(\right.
%       \sum_{n=m}^\infty\Cs_{n}^2\int_{\REAL}\big(\Hs_{n}^{(p)}\big)^2\expAW\ddv
%     }
%     - \sum_{\ell=p}^{m-1}\Cs_{\ell}^2\int_{\REAL}\big(\Hs_{\ell}^{(p)}\big)^2\expAW\ddv
%   \right)
%   \nonumber\\[0.5em]
%   &=
%   2^{m-p}\left
%   \sum_{n=p}^\infty\Cs_{n}^2\int_{\REAL}\big(\Hs_{n}^{(p)}\big)^2\expAW\ddv -
%   2^{m-p}\sum_{\ell=p}^{m-1}\Cs_{\ell}^2\int_{\REAL}\big(\Hs_{\ell}^{(p)}\big)^2\expAW\ddv
%   \nonumber\\[0.5em]
%   &=
%   2^{m-p}\int_{\REAL}\big(\varphi^{(p)}\big)^2\expAW\ddv -
%   2^{m-p}\sum_{\ell=p}^{m-1}\Cs_{\ell}^2\int_{\REAL}\big(\Hs_{\ell}^{(p)}\big)^2\expAW\ddv
%   \label{eq:2506b}
% \end{align}
% Finally, we reverse the inequality and estimate the additional terms
% to find that
\begin{align}
  \int_{\REAL}\big(\varphi^{(p)}\big)^2\expAW\ddv
  &\leq
  \frac{1}{2^{m-p}\,(m-p)!}\int_{\REAL}\big(\varphi^{(m)}\big)^2\expAW\ddv
  + \sum_{\ell=p}^{m-1}\Cs_{\ell}^2\int_{\REAL}\big(\Hs_{\ell}^{(p)}\big)^2\expAW\ddv
  \nonumber\\[0.5em]
  &=
  \frac{1}{2^{m-p}\,(m-p)!}\int_{\REAL}\big(\varphi^{(m)}\big)^2\expAW\ddv
  + \sum_{\ell=p}^{m-1}\Cs_{\ell}^2\,2^{p}\,\frac{\ell!}{(\ell-p)!}\int_{\REAL}\Hs_{\ell}^{2}\expAW\ddv
  \nonumber\\[0.5em]
  &=
  \frac{1}{2^{m-p}\,(m-p)!}\int_{\REAL}\big(\varphi^{(m)}\big)^2\expAW\ddv
  + 2^{p}\,\sqrt{\pi}\sum_{\ell=p}^{m-1}\,2^{\ell}\,\frac{(\ell!)^2}{(\ell-p)!}\,\Cs_{\ell}^2.
  \label{eq:2503a}
\end{align}
In particular, if $\varphi$ belongs to the space of polynomials of
degree at most $N$, we have $2n\leq 2N$, so that the relations
in~\eqref{eq:2501} can be adjusted to obtain the so called {\sl
  inverse inequality}
\begin{align}
 \int_{\REAL}\big(\varphi^{\prime}\big)^2\expAW\ddv\leq 2N
    \int_{\REAL}\varphi^2\expAW\ddv.
  \label{eq:2502a}
\end{align}
%% that now does not require restrictions on $C_0$.

\smallskip
Another useful inequality can be derived as follows.
First of all, from~\eqref{eq:0hermite:2} and~\eqref{eq:0hermite:3}, we know that:
\begin{align}
  2\vs\Hs_{n} 
  = \Hs_{n}^{\prime} + \Hs_{n+1}
  = 2n\Hs_{n-1} + \Hs_{n+1}
  \quad\forall n\geq 1.
  % \label{eq:0hermite:2}
\end{align}
Afterwords, we start by showing that:
\begin{align}
  \int_{\REAL}\vs^2\Hs_{n}^2\expAW\ddv 
  &= \int_{\REAL}n^2\Hs_{n-1}^2\expAW\ddv
  +\frac{1}{4}\int_{\REAL}\Hs_{n+1}^2\expAW\ddv\nonumber\\[0.5em]
  &= \sqrt{\pi}
  \left[
    n^2\,2^{n-1}\,(n-1)! + \frac{1}{4}\,2^{n+1}\,(n+1)!
  \right] = \sqrt{\pi}\left[
    2^{n-1}\,n\,n! + 2^{n-1}(n+1)n!
  \right]\nonumber\\[0.5em]
  &= \sqrt{\pi} 2^{n-1}(2n+1)n!
 \ \leq \ \sqrt{\pi}\,\frac{3}{4}\,2^{n+1}\,n\,n!
  = \frac{3}{4}\int_{\REAL}\big(\Hsp_{n}\big)^2\expAW\ddv,
  \qquad\forall n\geq 1,
\end{align}
where we noted that $2n+1\leq 2n+n=3n$, since $n\geq1$. 
The last equality follows from~\eqref{eq:005}. 
In short, we can write:
\begin{align}
  \int_{\REAL}\vs^2\Hs_{n}^2\expAW\ddv
  \leq \frac{3}{4}\int_{\REAL}\big(\Hsp_{n}\big)^2\expAW\ddv,
  \qquad\forall n\geq 1.
\end{align}

In general, let us suppose that $\varphi$ is a polynomial of degree
$N$ with $\Cs_{0}=0$.  
Thus, $\varphi$ has an expansion of the type
$\varphi=\sum_{n=1}^N\Cs_{n}\Hs_{n}$. 
For a given set of values $\alpha_n$, the following relation is a
consequence of the Schwartz inequality:
\begin{align}
  \left(\sum_{n=1}^{N}\alpha_n \right)^2 = \left(\sum_{n=1}^{N} 1\cdot \alpha_n \right)^2 \leq
   \sum_{n=1}^{N}1^2 \ \sum_{n=1}^{N}\alpha_n^2 = N \sum_{n=1}^{N}\alpha_n^2.
\end{align}
With the help of the above inequality, the orthogonality of the Hermite
polynomials implies that:
\begin{align}
  \int_{\REAL}\vs^2\varphi^2\expAW\ddv
  &= \int_{\REAL}\vs^2\Big(\sum_{n=1}^{N}\Cs_n\Hs_{n}\Big)^2\expAW\ddv
  \leq N\sum_{n=1}^{N}\Cs_{n}^2\int_{\REAL}\vs^2\Hs^2_{n}\expAW\ddv\nonumber\\[0.5em]
  &\leq \frac{3}{4}N\sum_{n=1}^{N}\Cs_{n}^2\int_{\REAL}\big(\Hsp_{n}\big)^2\expAW\ddv%\nonumber\\[0.5em]
  = \frac{3}{4}N\int_{\REAL}\big(\varphi'\big)^2\expAW\ddv,
  %%\qquad\textrm{~for~every~}\varphi\textrm{~with~}\Cs_{0}=0.
  \label{eq:2600}
\end{align}
which holds for every polynomial $\varphi$ with degree less or equal
to $N$ and $\Cs_{0}=0$.

\medskip
We end this preliminary section by introducing a few definitions
concerning the \emph{Hermite functions}, i.e., those functions that can be
written as a linear combination (finite or infinite) of the elements
of the \emph{Hermite basis functions} $\{\psi_{n}\}$.
Following the current literature, we will adopt a suitable notation in
order to distinguish the so-called \emph{symmetrically-weighted} (SW)
case, from the \emph{asymmetrically-weighted} (AW) one.
The reason of this setting will be made clear as we proceed with the
exposition.
We then consider the following definition:
%$\psi=\sum_{n=0}^\infty\Cs_{n}\psi_n$, where
\begin{align}
  \psi_n(\vs) 
  = \begin{cases}
    \gamma_{n}^{SW}\Hs_{n}(\vs)e^{-\vs^2\slash{2}} & \textrm{symmetrically-weighted~case},\\
    \gamma_{n}^{AW}\Hs_{n}(\vs)e^{-\vs^2}         & \textrm{asymmetrically-weighted~case},
  \end{cases}
  \label{eq:0245}
\end{align}
for some suitable choice of the real scalar coefficients
$\gamma_{n}^{SW}$ and $\gamma_{n}^{AW}$ (see below).
%%
%The Fourier coefficients $\Cs_{n}$ are obtained from the relation:
%\begin{align}
%  \Cs_{n} = \int_{\REAL}\psi\psi^n\ddv,
%\end{align}
%where $\psi^{n}$ for every integer $1\leq n\leq\infty$ belongs to the
%set of the dual basis functions defined by:
Besides, we introduce  the dual basis functions defined by:
\begin{align}
  \psi^n(\vs) 
  = \begin{cases}
    \widetilde{\gamma}_{n}^{SW}\Hs_{n}(\vs)e^{-\vs^2\slash{2}} & \textrm{symmetrically-weighted~case},\\
    \widetilde{\gamma}_{n}^{AW}\Hs_{n}(\vs)                  & \textrm{asymmetrically-weighted~case}.  
  \end{cases}
  \label{eq:0250}
\end{align}
The coefficients $\widetilde{\gamma}_n^{SW}$ and
$\widetilde{\gamma}_n^{AW}$ are obtained from the
orthogonality relation:
\begin{align}
  \left< \psi_n, \psi^m \right> = \delta_{n,m}.
  \label{eq:0255}
\end{align}
% To ease the exposition, we will remove the superscripts $AW$ and $SW$
% from the $\gamma$-coefficients.
%% 
We have:  
\begin{align}
  \gamma_{n}^{SW} =  \widetilde{\gamma}_{n}^{SW} = (\sqrt{\pi}2^n\,n!)^{-\frac{1}{2}},
  \label{eq:gamma_SW}
\end{align}
and
\begin{align}
  \gamma_{n}^{AW} = \,(\pi 2^n\,n!)^{-\frac{1}{2}},\qquad
  \widetilde{\gamma}_{n}^{AW} = (2^n\,n!)^{-\frac{1}{2}}.
  \label{eq:gamma_AW}
\end{align}

\section{Diffusive operators in the AW case}
\label{sec:diffusive:operators:AW}
Throughout the paper we will use indifferently the notation
$\partial\fs\slash{\partial\vs}$ and $\fs'$ to denote the partial
derivative of functions like $\fs(\vs)$ or $\fs(t,\vs)$, regardless of
their possible dependence on time.

We begin with the study of the second-order ($k=1$) differential
operator that appears in the Vlasov equation~\eqref{eq:Vlasov:1} and
the simplified model equation~\eqref{eq:model:1}.
In the asymmetric case, this operator can be decomposed as the
functional product of the two first-order operators:
\begin{align}
  \Ls = \frac{1}{2}\frac{\partial}{\partial\vs} + \vs\calI,\qquad\qquad
  \Lst = \frac{\partial}{\partial\vs},
  \label{eq:0000}
\end{align}
with $\calI$ the identity operator.
The second operator, i.e., $\Lst$, is just the derivative with respect
to the variable $\vs$.

\smallskip
We investigate the action of $\Lst\Ls$ on Hermite functions that we
write in the form:
\begin{align}
  \fs(\vs) = \hs(\vs)e^{-\vs^2},
  \label{eq:0005}
\end{align}
where $\hs$ is a generic polynomial.
%%
% To ease the exposition, we use $\hsp$, $\hspp$, $\cdots$, $\hs^{(k)}$
% to denote the first, second, \ldots, $k$-th derivative of $\hs$.
% %%
% Incidentally, $\hs$ can also be a function of time $t$ and position
% $x$, as it will be examined later on.
%%
For the operator $\Ls$, we have:
\begin{align}
  \Ls\fs 
  = \left(\frac{1}{2}\frac{\partial}{\partial\vs} + \vs\calI\right)\fs         
  %= \left(\frac{1}{2}\frac{\partial}{\partial\vs} + \vs\calI\right)\hs e^{-\vs^2} 
  = \frac{1}{2}\hsp e^{-\vs^2} - v\hs e^{-\vs^2} +  \vs\hs e^{-\vs^2}      
  = \frac{1}{2}\hsp e^{-\vs^2}.
  \label{eq:0010}
\end{align}
Clearly, $\Ls\fs$ is identically zero if $\hs$ is a constant.
Therefore, by taking $\hs=1$ we find that $\Ls(\expAW)=0$.

\smallskip
Similarly, for $k=2$ we have 
\begin{align}
  \Ls^2\fs
  = \Ls( \Ls\fs ) 
  = \Ls\left( \frac{1}{2}\hsp e^{-\vs^2} \right)
  = \frac{1}{4}\hspp e^{-\vs^2} - \frac{1}{2}v\hsp e^{-\vs^2} + \frac{1}{2}\vs\hsp e^{-\vs^2}
  = \frac{1}{4}\hspp e^{-\vs^2},
  \label{eq:0015}
\end{align}
and, in general, for $k\geq 2$ we have
\begin{align}
  \Ls^{k}\fs
  = \Ls( \Ls^{k-1}\fs ) 
  = \frac{1}{2^{k}}\hs^{(k)} e^{-\vs^2}.
  \label{eq:0025}
\end{align}
Equation~\eqref{eq:0025} can be proved recursively by
using~\eqref{eq:0010} for the first step, assuming that
$\Ls^{k-1}=(1\slash{2^{k-1}})\hs^{(k-1)}\expAW$ and applying the
definition of $\Ls$ given in~\eqref{eq:0000} to derive the relation at
step $k$.

\medskip
The combination of $\Ls$ and $\Lst$ provides the so called
second-order \emph{Lenard-Bernstein-like
  operator}~\cite{Lenard-Bernstein:1958}:
\begin{align}
  \Lst\Ls\fs 
  = \Lst\left(\frac{1}{2}\frac{\partial}{\partial\vs} + \vs\calI\right)\fs
  = \Lst\left( \frac{1}{2}\hsp e^{-\vs^2} \right)
  = \frac{1}{2}\hspp e^{-\vs^2} - \hsp \vs e^{-\vs^2}.
  \label{eq:0030}
\end{align}
Within the space of polynomials, $\Lst\Ls\fs$ is zero if and only if
$\hs$ is constant.
The combined operator is diffusive.
To prove this statement, we consider the time dependent problem for
the unknown function $\fs(\vs,t)=\hs(\vs,t)\expAW$:
\begin{align}
  \frac{\partial\fs}{\partial t} - \Lst\Ls\fs
  = \frac{\partial\fs}{\partial t} - \frac{\partial\Ls\fs}{\partial\vs}
  = 0,
  \label{eq:0035}
\end{align}
where again we assume that $\hs$ is a polynomial with respect to
$\vs$.
We multiply~\eqref{eq:0035} by $\hs$, integrate over $\REAL$, and,
then, integrate by parts the second integrand.
The boundary terms are zero since they can be expressed as a
polynomial multiplied by $e^{-\vs^2}$, which tends to zero for
$\ABS{\vs}\to\infty$.
Considering the expression of $\Ls\fs$ given in~\eqref{eq:0010}, we
obtain:
\begin{align}
  0
  &= \int_{\REAL}\left(\frac{\partial\fs}{\partial t} - \Lst\Ls\fs\right)\hs\ddv 
  = \int_{\REAL}\left(\frac{\partial\fs}{\partial t} - \frac{\partial\Ls\fs}{\partial\vs}\right)\hs\ddv 
  = \int_{\REAL}\frac{\partial\fs}{\partial t}\hs + \int_{\REAL}\big(\Ls\fs\big)\,\hsp\ddv
  -\left[\big(\Ls\fs\big)\hs\right]_{-\infty}^{+\infty}
  \nonumber\\[0.5em]
  &= \frac{1}{2}\frac{d}{\dt}\int_{\REAL}\hs^2e^{-\vs^2}\ddv + \frac{1}{2}\int_{\REAL}\big(\hsp\big)^2e^{-\vs^2}\ddv.
  \label{eq:0040}
\end{align}
From the equation above it follows that:
\begin{align}
  \frac{d}{\dt}\int_{\REAL}\hs^2 e^{-\vs^2}\ddv
  = -\int_{\REAL}\big(\hsp\big)^2e^{-\vs^2}\ddv \leq 0, 
  \label{eq:0045}
\end{align}
% which implies that~\eqref{eq:0035} is a parabolic equation in the weighted norm
% \begin{align*}
%   \norm{\fs}{\LTWO(\REAL)}^2 := \int_{\REAL}\abs{\fs}^2e^{+\vs^2}\ddv,
% \end{align*}
so that $\Lst\Ls\fs$ can be considered a dissipative operator for the
weighted $L^2(\REAL )$ norm.

\smallskip
Next, we repeat the same analysis for the fourth-order operator ($k=2$).
Consider again $\fs(\vs,t)=\hs(\vs,t)e^{-\vs^2}$ with $h$ polynomial,
and the time dependent problem:
\begin{align}
  \frac{\partial\fs}{\partial t} + \Lst^2\Ls^2\fs 
  = \frac{\partial\fs}{\partial t} + \frac{\partial^2\Ls^2\fs}{\partial\vs^2} = 0
  \label{eq:0050}
\end{align}
(note the change of sign with respect to Eq.~(\ref{eq:0035})).
As before, we multiply~\eqref{eq:0050} by $\hs$ and integrate over $\REAL$. 
Using the integration by parts (twice), we note that all the boundary
terms are zero since they always consist of a polynomial function in
$\vs$ multiplied by the Gaussian function $\expAW$, which tends to
zero for $\ABS{\vs}\to\infty$.
Omitting the boundary terms and using~\eqref{eq:0015} in the next
calculation, we obtain:
\begin{align}
  0 
  &= \int_{\REAL}\left( \frac{\partial\fs}{\partial t} + \frac{\partial^2\Ls^2\fs}{\partial\vs^2} \right)\hs\ddv
  % = \frac{1}{2}\frac{d}{\dt}\int_{\REAL}\hs^2 e^{-\vs^2}\ddv - \int_{\REAL}\frac{\partial\Ls^2\fs}{\partial\vs}\hsp\ddv\nonumber\\[0.5em]
  % &
  = \frac{1}{2}\frac{d}{\dt}\int_{\REAL}\hs^2e^{-\vs^2}\ddv + \int_{\REAL}\big(\Ls^2\fs\big)\,\hspp\ddv                \nonumber\\[0.5em]
  &= \frac{1}{2}\frac{d}{\dt}\int_{\REAL}\hs^2e^{-\vs^2}\ddv + \frac{1}{4}\int_{\REAL}\big(\hspp\big)^2e^{-\vs^2}\ddv.
  \label{eq:0055}
\end{align}
The equations above imply that $-\Lst^2\Ls^2\fs$ plays the role of a
diffusive term, since:
\begin{align}
  \frac{d}{\dt}\int_{\REAL}\hs^2e^{-\vs^2}\ddv =
  -\frac{1}{2}\int_{\REAL}\big(\hspp\big)^2e^{-\vs^2}\ddv\leq 0.
  \label{eq:0060}
\end{align}

The general case can be handled in a very similar way.
We write the time-dependent problem with the $2k$-th order operator as
follows:
\begin{align}
  \frac{\partial\fs}{\partial t} + (-1)^{k}\Lst^k\Ls^k\fs 
  = \frac{\partial\fs}{\partial t} + (-1)^{k}\frac{\partial^k\Ls^k\fs}{\partial\vs^k} 
  = 0\qquad(k\geq 1).
  \label{eq:0065}
\end{align}
Repeating the same arguments it follows that $-(-1)^k \Lst^k\Ls^k\fs$
is a diffusive operator.
Indeed, applying the integration by parts ($k$ times) and
recalling~\eqref{eq:0025}, yields:
\begin{align}
  -(-1)^{k}\int_{\REAL}\big(\Lst^k\Ls^k\fs\big)\hs\ddv
  &= -(-1)^{k}\int_{\REAL}\frac{\partial^k\Ls^k\fs}{\partial\vs^k}\,\hs\ddv
	\nonumber\\[0.5em]
  &=- (-1)^{k}(-1)^{k}\int_{\REAL}\big(\Ls^k\fs\big)\,\hs^{(k)}\ddv
  = -\frac{1}{2^{k}}\int_{\REAL}\big(\hs^{(k)}\big)^2\expAW\ddv
  \leq 0,
  \label{eq:0067}
\end{align}
where $\hs^{(k)}$ is the $k$-th derivative of $\hs$ with respect to
$\vs$.
The operators of order $2k$ for $k\geq 1$ so far examined are not
strictly negative definite, since their kernel is not empty.

\smallskip

% The operators of order $2k$ for $k\geq 1$ so far examined are not
% strictly negative-definite, since their kernel is not empty.
% %% 
% They have been constructed with this property for the purpose of
% stabilizing certain equations and preserving at the same time the
% lower modes.\textcolor{red}{Queste ultime due frasi sono un po' 'out
%   of the blue' e le toglierei. Soprattutto l'ultima visto che non
%   considera le leggi di conservazione e la filamentazione.}

\section{Action of the diffusive operators in the AW Hermite case}
\label{sec:action:operators:AW}
Consider~\eqref{eq:0030} in terms of the Hermite functions' basis.
A direct calculation yields:
\begin{align}
  \Lst\Ls\psi_{n} 
  &= \Lst\Ls\big( \gamma^{AW}_{n}\Hs_{n}\expAW \big) 
  = \frac{\gamma^{AW}_{n}}{2}\big( \Hs_{n}^{\pp} - 2\vs\Hs_{n}^{\prime} \big) \expAW
  = \frac{\gamma^{AW}_{n}}{2}\big( -2n\Hs_{n}\expAW \big)
	\nonumber\\[0.5em]
  &= -n\gamma^{AW}_{n}\Hs_{n}\expAW
  = -n\psi_{n},	
  \label{eq:0030a}
\end{align}
where we used the differential equation~\eqref{eq:0075}.
In other words, the function $\psi_{n}$ is the eigenfunction of the
differential operator $\Lst\Ls$ with eigenvalue $-n$.
As the corresponding eigenvalue is zero for $n=0$, it follows that
$\Lst\Ls$ acts on Hermite functions without altering the equation for the first
Hermite coefficient $\Cs_{0}$.
This is a further confirmation of the diffusive nature of the operator
regarding the Hermite modes that are higher than $1$.

\smallskip
%%
%%We have demonstrated that the element of the AW Hermite functions' basis
%%$\psi_n$ is an eigenfunction of the partial differential operator
%%$\Lst\Ls$ with eigenvalue $-n$.
%%
A similar relation holds also for $\Lst^2\Ls^2$ and for
the more general operator $\Lst^{k}\Ls^{k}$.
First, we consider the case $k=2$.
Using~\eqref{eq:0015} with $\hs=\Hs_{n}$, a straightforward calculation
yields:
\begin{align}
  \Lst^2\Ls^2\psi_{n} 
  = \Lst^2\Ls^2\big( \gamma^{AW}_{n}\Hs_{n}\expAW \big)
  = \Lst^{2}\Big( \gamma^{AW}_{n} \frac{1}{4}\Hs_{n}^{\pp}\expAW \Big)
  = \frac{\gamma^{AW}_{n}}{4}\Big( \Hs_{n}^{\pp}\expAW \Big)^{\pp}.
  \label{eq:0090b}
\end{align}
To compute the last term in the equation above, we proceed in two
steps, starting from the first derivative of $\Hs^{\pp}\expAW$.
Using~\eqref{eq:0075}, we have that:
% \begin{align}
%   \left( \Hs_{n}^{\pp}\expAW \right)^{\prime}
%   &= \left( \big( 2\vs\Hs_{n}^{\prime} -2n\Hs_{n} \big) \expAW \right)^{\prime}
% 	 \nonumber\\[0.5em]
%   &= \left( 2\Hs_{n}^{\prime} + 2\vs\Hs_{n}^{\pp} - 2n\Hs_{n}^{\prime} \right) \expAW
%   -2\vs \left( 2\vs\Hs_{n}^{\prime} - 2n\Hs_{n} \right) \expAW                       
%   \nonumber\\[0.5em]
%   &= \left( 2\Hs_{n}^{\prime} + 2\vs(2\vs\Hs_{n}^{\prime} - 2n\Hs_{n}) - 2n\Hs_{n}^{\prime} 
%     -4\vs^2\Hs_{n}^{\prime} + 4n\vs\Hs_{n} \right) \expAW                         
%   = 2(1-n)\Hs_{n}^{\prime} \expAW.
%   \label{eq:0100}
% \end{align}
\begin{align}
  \left( \Hs_{n}^{\pp}\expAW \right)^{\prime}
  &= \left( \big( 2\vs\Hs_{n}^{\prime} -2n\Hs_{n} \big) \expAW \right)^{\prime}
	 \nonumber\\[0.5em]
  &= \left( 2\Hs_{n}^{\prime} + 2\vs\Hs_{n}^{\pp} - 2n\Hs_{n}^{\prime} \right) \expAW
  -2\vs \left( 2\vs\Hs_{n}^{\prime} - 2n\Hs_{n} \right) \expAW                       
  \nonumber\\[0.5em]
  &= \left( 2\Hs_{n}^{\prime} + 2\vs\Hs_{n}^{\pp} - 2n\Hs_{n}^{\prime} \right) \expAW
  -2\vs\Hs_{n}^{\pp} \expAW                       
  \nonumber\\[0.5em]
  %% &= \left( 2\Hs_{n}^{\prime} + 2\vs(2\vs\Hs_{n}^{\prime} - 2n\Hs_{n}) - 2n\Hs_{n}^{\prime} 
  %%   -4\vs^2\Hs_{n}^{\prime} + 4n\vs\Hs_{n} \right) \expAW                         
  & = 2(1-n)\Hs_{n}^{\prime} \expAW.
  \label{eq:0100}
\end{align}
Using again~\eqref{eq:0075}, we have that:
\begin{align}
  \left( \Hs_{n}^{\prime}e^{-\vs^2} \right)^{\prime}
  = \Hs_{n}^{\pp}e^{-\vs^2} - 2\vs\Hs_{n}^{\prime}e^{-\vs^2}
  = \left( \Hs_{n}^{\pp} -2\vs\Hs_{n}^{\prime} \right)e^{-\vs^2}
  = -2n\Hs_{n}e^{-\vs^2}.
  \label{eq:0090}
\end{align}
Hence, the second derivative of $\Hs^{\pp}\expAW$ with respect to
$\vs$ is readily given by collecting the results of~\eqref{eq:0100}
and~\eqref{eq:0090}, and reads as:
\begin{align}
  \left( \Hs_n^{\pp}\expAW \right)^{\pp} 
  &= \left( \left( \Hs_{n}^{\pp}\expAW \right)^{\prime} \right)^{\prime}
  = \left( 2(1-n)\Hs_{n}^{\prime}\expAW \right)^{\prime}
  %%&= 2(1-n)\left( \Hs_{n}^{\prime}\expAW \right)^{\prime}
  = 4n(n-1)\Hs_{n}\expAW.
  \label{eq:0105}
\end{align}
Replacing~\eqref{eq:0105} in~\eqref{eq:0090b}, finally yields:
\begin{align}
  \Lst^2\Ls^2\psi_{n} 
  = \frac{\gamma^{AW}_{n}}{4}\,4n(n-1)\Hs_{n}\expAW
  = n(n-1)\,\gamma^{AW}_{n}\Hs_{n}\expAW
  = n(n-1)\psi_{n},
\end{align}
which shows that $\psi_n$ is an eigenfunction of $\Lst^2\Ls^2$
corresponding to the eigenvalue $n(n-1)$.
Note that such eigenvalue is zero for $n=0$ and $n=1$, which means
that the fourth-order operator $\Lst^2\Ls^2$
does not modify the equations for the first two modes of the AW Hermite expansion of
$\fs$. 

\medskip
Repeating the same argument for a general integer $k\geq1$, we find
out that:
\begin{align}
  \Lst^k\Ls^k\psi_{n} 
  = (-1)^{k}\,n(n-1)\ldots(n-(k-1))\psi_{n} 
  = (-1)^k\frac{n!}{(n-k)!}\psi_{n}.
  \label{eq:0076}
\end{align}
Therefore, we conclude that every element of the AW Hermite function's
basis is an eigenfunction of the $2k$-th operator $\Lst^{k}\Ls^{k}$
with eigenvalue $(-1)^{k}\,n!\slash{(n-k)!}$, which takes the value of
zero for $0\leq n\leq k-1$.
%%
% We can similarly conclude that the $2k$-th operator does not modify
% the first $k$ modes of the AW Hermite expansion of $\fs$.

\medskip
We conclude this section by investigating the action of the
Lenard-Bernstein-like operators on Hermite functions expressed as
linear combinations of the AW Hermite functions' basis and the
implications on the conservation properties of the discretization.
Similar topics were considered in the more specific context of
Vlasov-based models,
cf.~\cite{Delzanno:2015,Camporeale-Delzanno-Bergen-Moulton:2015}.

To this end, we consider again the expansion $\fs(\vs)=h(\vs)\expAW$,
where the polynomial function is given by (see~\eqref{eq:0052}):
\begin{align}
  \hs = \sum_{n=0}^\infty\Cs_{n}\Hs_{n}.
  \label{eq:0070}
\end{align}
By multiplying and dividing by the normalization factor
$\gamma^{AW}_n$, and, then, using the definition of the AW basis
(see~\eqref{eq:0245}-\eqref{eq:0250}) we find that:
\begin{align}
  \fs 
  = \hs\expAW 
  = \left(\sum_{n=0}^\infty\Cs_{n}\Hs_{n}\right)\expAW
  = \sum_{n=0}^\infty\frac{\Cs_{n}}{\gamma^{AW}_n}\left(\gamma^{AW}_n\Hs_{n}\expAW\right)
  = \sum_{n=0}^\infty\Csb_{n}\psi_{n},
  \label{eq:0070a}
\end{align}
where $\Csb_{n}=\Cs_{n}\slash{\gamma^{AW}_{n}}$.
Since $\psi_{n}$ is an eigenfunction of the generalized Lenard-Bernstein
operators, we obtain the following relations:
\begin{align}
  \Lst\Ls\fs 
  &= \sum_{n=0}^{\infty}\Csb_{n}\Lst\Ls\psi_{n}
  = \sum_{n=0}^{\infty}(-n)\Csb_{n}\psi_{n},\\[0.5em]
  %% --------------------------------------------
  \Lst^2\Ls^2\fs^2 
  &= \sum_{n=0}^{\infty}\Csb_{n}\Lst^2\Ls^2\psi_{n}
  = \sum_{n=0}^{\infty}n(n-1)\Csb_{n}\psi_{n},\\[0.5em]
  %% --------------------------------------------
  \ldots\nonumber\\[0.5em]
  %% --------------------------------------------
  \Lst^k\Ls^k\fs 
  &= \sum_{n=0}^{\infty}\Csb_{n}\Lst^k\Ls^k\psi_{n}
  = \sum_{n=0}^{\infty} (-1)^{k}\frac{n!}{(n-k)!}\Csb_{n}\psi_{n}.
\end{align}
From the identities above, it follows immediately that:
\begin{align}
  \Lst\Ls\fs &= \sum_{n=0}^\infty\Ds^{(1)}_{n}\psi_{n}
  \qquad\textrm{with~}\Ds^{(1)}_{n} = -n\Csb_{n},\\[0.5em]
  %% -----------------------------------------------------------
  \Lst^2\Ls^2\fs &= \sum_{n=0}^\infty\Ds^{(2)}_{n}\psi_{n}
  \qquad\textrm{with~}\Ds^{(2)}_{n} = n(n-1)\Csb_{n},\\[0.5em]
  %% -----------------------------------------------------------
  \ldots\nonumber\\[0.5em]
  %% -----------------------------------------------------------
  \Lst^{k}\Ls^{k}\fs &= \sum_{n=0}^\infty\Ds^{(k)}_{n}\psi_{n}
  \qquad\textrm{with~}\Ds^{(k)}_{n} = (-1)^{k}\frac{n!}{(n-k)!}\Csb_{n}.
\end{align}
By definition, it holds that
$\Ds^{(k)}_{0}=\Ds^{(k)}_{1}=\ldots=\Ds^{(k)}_{k-1}=0$ for a generic
$k\geq1$. The case $k=3$ corresponds to the operator used in Refs. \cite{Camporeale-Delzanno-Bergen-Moulton:2016,Delzanno:2015}.
%% , we can start the above summations from $n=0$.

% Consistently with the previous comments, we can set $\Ds^{(1)}_{0}=0$,
% for $k=1$; $\Ds^{(2)}_{0}=\Ds^{(2)}_{1}=0$, for $k=2$.

%% \subsection{Conservation properties}
%%
\medskip
Using the properties that we have established so far, we are able to
prove some conservation properties for problems of parabolic type like
those considered in~\eqref{eq:0035} (using $\Lst\Ls\fs$),
~\eqref{eq:0050} (using $-\Lst^2\Ls^2\fs$), ~\eqref{eq:0065} (using
$-(-1)^k\Lst^k\Ls^k\fs$).
%%
%% Going back to the model problems of 
%% examined in the previous section,
%%
The \emph{mass conservation} for a distribution function $\fs(t,\vs)$
is expressed by:
\begin{align}
  \frac{d}{\dt}\int_{\REAL}\fs\ddv = 0.
  \label{eq:0125}
\end{align}
%%
% Then, we consider the three cases in which $\fs$ is the solution of
% ~\eqref{eq:0035} (using $\Lst\Ls\fs$), ~\eqref{eq:0050} (using
% $-\Lst^2\Ls^2\fs$), ~\eqref{eq:0065} (using $-(-1)^k\Lst^k\Ls^k\fs$).
%%
%\SubItem{Mass conservation for Eq.~\eqref{eq:0035} -- using diffusion
%  operator $\Lst\Ls\fs$}
\smallskip
In the first case, we integrate~\eqref{eq:0125} on $\REAL$,
use~\eqref{eq:0035}, apply the fundamental theorem of calculus and
substitute the expression of $\Ls\fs$ in~\eqref{eq:0010} to obtain:
\begin{align}
  \frac{d}{\dt}\int_{\REAL}\fs\ddv 
  = \int_{\REAL}\Lst\Ls\fs\ddv
  = \int_{\REAL}\frac{\partial\Ls\fs}{\partial\vs}\ddv
  = \left[\Ls\fs\right]^{\infty}_{-\infty}
  = \frac{1}{2}\left[ \hsp e^{-\vs^2} \right]^{\infty}_{-\infty}
  = 0,
  \label{eq:0130}
\end{align}
since $e^{-\vs^2}$ times a polynomial of any degree tends to zero for
$\vs\to\pm\infty$.

In the second case, we integrate~\eqref{eq:0125} on $\REAL$,
use~\eqref{eq:0050}, and apply the fundamental theorem of calculus to
obtain:
\begin{align}
  \frac{d}{\dt}\int_{\REAL}\fs\ddv 
  = -\int_{\REAL}\Lst^2\Ls^2\fs\ddv
  = -\int_{\REAL}\frac{\partial}{\partial\vs}(\Lst\Ls^2\fs )\ddv
  = -\left[ \Lst\Ls^2\fs \right]^{\infty}_{-\infty}.
  \label{eq:0135}
\end{align}
Furthermore, by using~\eqref{eq:0015}, we find that:
\begin{align}
  \Lst\Ls^2\fs
  = \frac{\partial\Ls^2\fs}{\partial\vs}
  = \frac{1}{4}\frac{\partial}{\partial\vs}\left( \hspp\expAW \right)
  = \frac{1}{4}\big( \hsppp - 2\vs\hspp \big) \expAW.
  \label{eq:0140}
\end{align}
Therefore, the last term above provides zero in~\eqref{eq:0135}, since
the Gaussian function $\expAW$ multiplied by any polynomial tends to
zero for $\vs\to\pm\infty$.

Finally, to obtain the general result for $\Lst^k\Ls^k\fs$, we
integrate~\eqref{eq:0125} on $\REAL$, use~\eqref{eq:0065}, and apply
the fundamental theorem of calculus.
We obtain:
\begin{align}
  \frac{d}{\dt}\int_{\REAL}\fs\ddv 
  &= -(-1)^k\int_{\REAL}\Lst^k\Ls^k\fs\ddv
  = -(-1)^k\int_{\REAL}\frac{\partial}{\partial\vs}( \Lst^{k-1}\Ls^k\fs )\ddv \nonumber\\[0.5em]
  &= -(-1)^k\left[ \Lst^{k-1}\Ls^k\fs \right]^{\infty}_{-\infty}
  = 0,
  \label{eq:0145}
\end{align}
since we can prove recursively that $\Lst^{k-1}\Ls^{k}\fs$ is equal to
a polynomial times $e^{-\vs^2}$, which tends to zero for
$\vs\to\pm\infty$.

\smallskip
Another important issue is the \emph{momentum conservation}, which is
expressed by:
\begin{align}
  \frac{d}{\dt}\int_{\REAL}\vs\fs\ddv = 0.
  \label{eq:0150}
\end{align}
We start by noting that there is no momentum conservation for the
operator $\Lst\Ls$.  
We then consider the two other cases in which $\fs$ is the solution
of~\eqref{eq:0050} (using$-\Lst^2\Ls^2\fs$), and~\eqref{eq:0065}
(using $-(-1)^k\Lst^k\Ls^k\fs$).

\smallskip
In the first case, momentum conservation is achieved because, in view
of~\eqref{eq:0050}, we know that:
\begin{align}
  \frac{d}{\dt}\int_{\REAL}\vs\fs\ddv 
  = -\int_{\REAL}\vs\frac{\partial^2\Ls^2\fs}{\partial\vs^2}\ddv.
  \label{eq:0155}
\end{align}
Then, we integrate by parts the right-hand side, apply the fundamental
theorem of calculus and arrive at:
\begin{align}
  \frac{d}{\dt}\int_{\REAL}\vs\fs\ddv 
  = \int_{\REAL}\frac{\partial\Ls^2\fs}{\partial\vs}\ddv
  -\left[ \vs\frac{\partial\Ls^2\fs}{\partial\vs} \right]^{\infty}_{-\infty}
  =\left[ \Ls^2\fs \right]^{\infty}_{-\infty} 
  -\left[ \vs\frac{\partial\Ls^2\fs}{\partial\vs} \right]^{\infty}_{-\infty}
  = 0,
  \label{eq:0155a}
\end{align}
As in the previous situations, the arguments in the square brackets
are of the form of a polynomial multiplied by the Gaussian function
$\expAW$.

\smallskip
Through very similar steps, we can easily arrive at a general
statement regarding the conservation of the $m$-th moment,
$m\geq 1$.
Indeed, we have:
\begin{align}
  \frac{d}{\dt}\int_{\REAL}\vs^m\fs\ddv = 0
  \label{eq:0151}
\end{align}
in presence of the operator $\Lst^k\Ls^k\fs$, and provided that the
condition $k>m$ is satisfied.
The conservation of the velocity moments of the distribution function
$\fs$ implies the conservation of physical quantities such as momentum
and energy in Vlasov models.
We will discuss this topic at the beginning of  
Section~\ref{sec:full-discretization:Vlasov-Poisson}.
 
\section{Diffusive operators in the SW case}
\label{sec:diffusive:operators:SW}
Differently from the AW case, the generalized Lenard-Bernstein operators that
we consider in the SW case read as follows:
\begin{align}
  \Ls  = \frac{\partial}{\partial\vs} + \vs\calI,\qquad
  \Lst = \frac{\partial}{\partial\vs} - \vs\calI .
  \label{eq:0160}
\end{align}
We investigate the action of the $\Lst\Ls$ operator on Hermite
functions of the form $\fs=\hs\expSW$, where $\hs$ is once again a
polynomial in $\vs$.
The weighted $L^2$ inner product for such functions is:
\begin{align} 
  \big(\fs,\gs\big) 
  = \int_{\REAL}\fs\gs\ddv
  = \int_{\REAL}\hs_{\fs}\hs_{\gs}\expAW\ddv,
  \label{eq:0165}
\end{align}
where $\fs=\hs_{\fs}\expSW$ and $\gs=\hs_{\gs}\expSW$, and $\hs_{\fs}$
and $\hs_{\gs}$ are polynomials. 
This somehow justifies the adoption of the term ``symmetric".

\smallskip

The results will be analogous to those presented in the previous
sections.
We briefly review the main points.
From straightforward calculations it follows that:
\begin{align}
  \Ls\fs 
  &= \fs^{\prime}+\vs\fs 
  = \hs'\expSW - \vs\hs\expSW + \vs\hs\expSW
  = \hsp\expSW,\label{eq:0165a}\\[0.5em]
  \Lst\Ls\fs 
  &= \Lst\big(\hsp\expSW\big) = \big(\hspp-2\vs\hsp\big)\expSW.\label{eq:0165b}
\end{align}
%%
% A similar calculation yields that
% \begin{align}
%   \Lst\Ls\Hs_{n}\expSW = \big(\Hspp_{n}-2\vs\Hsp_{n}\big)\expSW = -2n\Hs_{n}\expSW.
%   \label{eq:0165b}
% \end{align}

These relations imply that the operator $\Lst\Ls$ is diffusive.
In fact, consider again the time dependent problem:
\begin{align}
  \frac{\partial\fs}{\partial t} - \Lst\Ls\fs = 0,
  \label{eq:0170}
\end{align}
where, now, we choose $f(v,t)=h(v,t)\expSW$.
We multiply equation~\eqref{eq:0170} by $\fs$ and integrate over
$\REAL$. 
Thus, we end up with the equality:
\begin{align}
  \int_{\REAL}\left( \frac{\partial\fs}{\partial t} - \Lst\Ls\fs \right)\fs\ddv = 0,
  \label{eq:0175}
\end{align}
and using the definition of $\Lst$ given in~\eqref{eq:0160}, we have that
\begin{align}
  \frac{1}{2}\int_{\REAL}\frac{\partial}{\partial t}\big(\hs^2\expAW\big)\ddv
  - \int_{\REAL}\big( (\Ls\fs)^{\prime} - \vs\Ls\fs \big)\fs\ddv = 0,
  \label{eq:0180}
\end{align}
where again we denoted the derivative with respect to $\vs$ of
$\Ls\fs$ by $(\Ls\fs)^{\prime}$.
Then, we integrate by parts the second integral of~\eqref{eq:0180} and
note again that the boundary terms for $\vs\to\pm\infty$ are zero.
This leads us to:
\begin{align}
  0
  &= \frac{1}{2}\frac{d}{\dt}\int_{\REAL}\hs^2 e^{-\vs^2}\ddv
  + \int_{\REAL} \big(\Ls\fs\big)\,\fs^{\prime}\ddv
  - \left[ (\Ls\fs)\fs \right]^{\infty}_{-\infty}
  + \int_{\REAL}\vs\big(\Ls\fs\big)\,\fs\ddv 
  \nonumber\\[0.5em]
  &= \frac{1}{2}\frac{d}{\dt}\int_{\REAL}\hs^2 e^{-\vs^2}\ddv
  + \int_{\REAL}\big(\Ls\fs\big)\big( \fs^{\prime}+\vs\fs \big)\ddv
  = \frac{1}{2}\frac{d}{\dt}\int_{\REAL}\hs^2 e^{-\vs^2}\ddv
  + \int_{\REAL}\big( \Ls\fs \big)^2\ddv.
  \label{eq:0185}
\end{align}
The last relation shows that the operator $\Lst\Ls$ introduces a
dissipation.

\smallskip 
The same result holds for the fourth-order operator and the related
time dependent problem:
\begin{align}
  \frac{\partial\fs}{\partial t} + \Lst^2\Ls^2\fs = 0.
  \label{eq:0200}
\end{align}
Here, the proof is a bit more involved, but still elementary.
We first note that
$\Lst^2\gs=\Lst(\Lst\gs)=\Lst(\gs^{\prime}-\vs\gs)$, from which it
follows that:
\begin{align}
  \Lst^2\gs 
  = ( \gs^{\prime}-\vs\gs )^{\prime} - \vs( \gs^{\prime} - \vs\gs )
  = \gs^{\pp} -2\vs\gs^{\prime} + (\vs^2-1)\gs,
  \label{eq:0205}
\end{align}
and
\begin{align}
  \Ls^2\fs 
  = \fs^{\pp} + 2\vs\fs^{\prime} + ( \vs^2+1 )\fs
  = \fs^{\pp} + 2\big(\vs\fs\big)^{\prime} + ( \vs^2-1 )\fs.
  \label{eq:0210}
\end{align}

By multiplying equation~\eqref{eq:0200} by $\fs$ and integrating over
$\REAL$, we find that:
\begin{align}
  \frac{1}{2}\frac{d}{\dt}\int_{\REAL}\hs^{2}e^{-\vs^2}\ddv 
  + \int_{\REAL}\big(\Lst^2\Ls^2\fs\big)\fs\ddv = 0.
  \label{eq:0215}
\end{align}
From straightforward calculations using integration by parts,
\eqref{eq:0205} (with $\gs=\Ls^2\fs$) and~\eqref{eq:0210}, we get the following relation:
\begin{align}
  \int_{\REAL}\big(\Lst^2\Ls^2\fs\big)\fs\ddv 
  &= \int_{\REAL}\left(
    \big( \Ls^2\fs \big)^{\pp} - 2\vs\big( \Ls^2\fs \big)^{\prime} + (\vs^2-1)\Ls^2\fs
  \right)\fs\ddv                                           \nonumber\\[0.5em]
  &= \int_{\REAL}\big( \Ls^2\fs \big)\,\fs^{\prime\prime}\ddv
  + 2\int_{\REAL}\big(\Ls^2\fs\big)\,\big( \vs\fs \big)^{\prime}\ddv
  + \int_{\REAL} (\vs^2-1)\big(\Ls^2\fs\big)\,\fs\ddv        \nonumber\\[0.5em]
  &= \int_{\REAL} \big( \Ls^2\fs \big)\left(
    \fs^{\pp} + 2\big( \vs\fs \big)^{\prime} + (\vs^2-1)\fs
  \right)\ddv                                              
  = \int_{\REAL} \big( \Ls^2\fs \big)^2\ddv.
  \label{eq:0220}
\end{align}
Therefore, also this time-dependent equation is dissipative, from the
viewepoint of the $L^2(\REAL )$-weighted norm.

\smallskip
In general, we may consider the time dependent problem:
\begin{align}
  \frac{\partial\fs}{\partial t} + (-1)^{k}\Lst^{k}\Ls^{k}\fs = 0.
  \label{eq:0225}
\end{align}
With the same considerations as above, we find the relation:
\begin{align}
  \frac{d}{\dt}\int_{\REAL}\hs^{2}\expAW\ddv 
  = -\int_{\REAL}\big(\Ls^k\fs\big)^{2}\ddv \leq 0,
  \label{eq:0230}
\end{align}
which shows the dissipative nature of the second term
of~\eqref{eq:0225}.

\smallskip
Regarding the expansion in the Hermite basis functions, after the
application of the diffusive operators, we also obtain straightforward
results.
First, we write the function $\fs(\vs)=h(\vs)\expSW$ in the SW Hermite
basis, by using the expansion:
\begin{align}
  \hs = \sum_{n=0}^\infty\Cs_{n}\Hs_{n},
  \qquad\textrm{with}\quad
  \Cs_{n} = \gamma^{SW}_n\tilde{\gamma}^{SW}_{n}\int_{\REAL}\hs\Hs_{n}\expAW\ddv.
  \label{eq:0071}
\end{align}
The corresponding coefficients $\Ds_{m}^{(1)}$ are such that:
\begin{align}
  \Lst\Ls\fs 
  = \Lst\left(\sum_{n=0}^\infty\Cs_{n}\Hs_{n}^{\prime}e^{-\vs^2/2}\right)
  = \sum_{m=1}^\infty\Ds_{m}^{(1)}\Hs_{m}e^{-\vs^2/2},
  \label{eq:0085}
\end{align}
which allows us to express $\Lst\Ls\fs$ in terms of Hermite functions.
In practice:
\begin{align}
  \Lst\Ls\left(\Hs_{n}\expSW\right)
  &=\Lst\left(\Hs_{n}^{\prime}\expSW\right)
  = \Hs_{n}^{\pp}e^{-\vs^2\slash{2}} 
  - \vs\Hs_{n}^{\prime}\expSW
  - \vs\Hs_{n}^{\prime}\expSW %
  \nonumber\\[0.5em]
  &= \Big( \Hs_{n}^{\pp} - 2\vs\Hs_{n}^{\prime} \Big)\expSW 
  = -2n\Hs_{n}\expSW,
  \label{eq:0190}
\end{align}
where we used again the differential equation for Hermite polynomials
(see~\eqref{eq:0075}).
%$\Hs_{n}^{\pp}-2\vs\Hs_{n}^{\prime}+2n\Hs_{n}=0$.
%%
Therefore,  we obtain:
\begin{align}
  \Ds_{m}^{(1)} = -2m\Cs_m
  \qquad\textrm{for~~}m\geq 1,
  \label{eq:0195}
\end{align}

%\NewItem{Expansion of $\Lst^k\Ls^k\fs$ using Hermite functions.}

Going to the general case, we want to compute the coefficients
$\Ds_{m}^{(k)}$ such that:
\begin{align}
  \Lst^k\Ls^k\fs 
  = \sum_{m=0}^{\infty}\Ds_{m}^{(k)}\Hs_{m}e^{-\vs^2/2},
  \label{eq:0098}
\end{align}
One finally obtains:
\begin{align}
  \Ds_{m}^{(k)} = (-1)^{k}2^{k}m(m-1)\cdots(m-k)\Cs_m = (-1)^k 2^{k}\frac{m!}{(m-k-1)!}\Cs_m
  \qquad\textrm{for~~}m\geq 0.
  \label{eq:0122}
\end{align}
As for the AW case, the first $k+1$ coefficients are automatically
zero.
%%
%%%%%%%%%%%%%%%%%%%%%%%%%%%%%%%%%%%%%%%%%%%%%%%%%%%%%%%%%%%%%%%%%%%%%%%%%%%%%%%

\section{Action of the diffusive operators in the SW Hermite case}
\label{sec:action:operators:SW}

We recall that $\Ls=\partial\slash{\partial\vs}+\vs\calI$ and
$\Lst=\partial\slash{\partial\vs}-\vs\calI$.
Consider the SW Hermite basis functions:
$\psi_{n}=\gamma^{SW}_n\Hs_{n}\expSW$.
A straightforward calculation yields:
\begin{align}
  \Ls\psi_n 
  &= \gamma^{SW}_n\left(\frac{\partial}{\partial\vs} + \vs\calI\right)\Hs_{n}\expSW\nonumber\\[0.5em]
  &= \gamma^{SW}_n\left(\Hsp_{n}\expSW-\vs\Hs_{n}\expSW+\vs\Hs_{n}\expSW\right)
  = \gamma^{SW}_n\Hsp_{n}\expSW.
\end{align}
Using the result above  we
obtain:
\begin{align}
  \Ls^2\psi_n
  & = \Ls\big(\Ls\psi_n\big) = \Ls\left( \gamma^{SW}_n\Hsp_{n}\expSW \right)
  = \gamma^{SW}_n\left(\frac{\partial}{\partial\vs} + \vs\calI\right)\Hsp_{n}\expSW\nonumber\\[0.5em]
  &= \gamma^{SW}_n\left(\Hspp_{n}\expSW-\vs\Hsp_{n}\expSW+\vs\Hsp_{n}\expSW\right)
  = \gamma^{SW}_n\Hspp_{n}\expSW.
\end{align}
A simple recursive argument allows us to prove the formula for a
generic $k$:
\begin{align}
  \Ls^{k}\psi_n = \gamma^{SW}_n\Hs_{n}^{(k)}\expSW,
\end{align}
where we recall that $\Hs^{(k)}_n=d^{k}\Hs_{n}\slash{d\vs^k}$.
Indeed, we have already proved that the formula is true for $k=1$ and
$k=2$.
Since $\Ls^{k-1}\psi_n=\gamma^{SW}_n\Hs_{n}^{(k-1)}\expSW$, a
straightforward calculation yields:
\begin{align}
  \Ls^{k}\psi_n 
  & = \Ls\big(\Ls^{(k-1)}\psi_n\big) = \Ls\left( \gamma^{SW}_n\Hs^{(k-1)}_{n}\expSW \right)
  = \gamma^{SW}_n\left(\frac{\partial}{\partial\vs} + \vs\calI\right)(\Hs^{(k-1)}_{n}\expSW)\nonumber\\[0.5em]
  &= \gamma^{SW}_n\left(\Hs^{(k)}_{n}\expSW -\vs\Hs^{(k-1)}_{n}\expSW+\vs\Hs^{(k-1)}_{n}\expSW\right)
  = \gamma^{SW}_n\Hs^{(k)}_{n}\expSW.
\end{align}
Now, we compute the action of $\Lst$, $\Lst^2$, and $\Lst^{k}$ on
$\Ls\psi_n$, $\Ls^2\psi_n$, and $\Lst^{k}\psi_n$, respectively.
In the first case, we recover the relation:
\begin{align}
  \Lst\Ls\psi_n 
  &= \gamma^{SW}_n\left(\frac{\partial}{\partial\vs} - \vs\calI\right)(\Hsp_{n}\expSW)
  =  \gamma^{SW}_n\left(\Hspp_{n}\expSW -\vs\Hsp_{n}\expSW-\vs\Hsp_{n}\expSW\right)\nonumber\\[0.5em]
  &= \gamma^{SW}_n\left(\Hspp_{n}-2\vs\Hsp_{n}\right)\expSW
  =  \gamma^{SW}_n(-2n)\Hs_{n} \expSW = -2n\psi_{n}, \qquad n\geq 1.
\end{align}
In the second case, first we obtain: 
\begin{align}
  \Lst\Ls^2\psi_n 
  &= \gamma^{SW}_n\left( \frac{\partial}{\partial\vs} - \vs\calI\right)(\Hspp_{n}\expSW)
  =  \gamma^{SW}_n\left( \big( \Hspp_{n} \big)^{\prime}\expSW -\vs\Hspp_{n}\expSW -\vs\Hspp_{n}\expSW \right)\nonumber\\[0.5em]
  &= \gamma^{SW}_n\left( \big( \Hspp_{n} \big)^{\prime} -2\vs\Hspp_{n} \right)\expSW
  =  \gamma^{SW}_n\left( \big(2\vs\Hsp_{n}-2n\Hs_{n}\big)^{\prime} -2\vs\Hspp_{n} \right)\expSW\nonumber\\[0.5em]
  &= \gamma^{SW}_n\left( 2\Hsp_{n}+2\vs\Hspp_{n}-2n\Hsp_{n}-2\vs\Hspp_{n} \right)\expSW
  =  \gamma^{SW}_n\,2(1-n)\,\Hsp_{n}\expSW , \quad n\geq 2,
\end{align}
and then:
\begin{align}
  \Lst^2\Ls^2\psi_n 
  &= \Lst(\Lst\Ls^2\psi_n) 
  = \Lst\big(\gamma^{SW}_n\, 2(1-n)\Hsp_n\expSW\big)\nonumber\\[0.5em]
  &= \gamma^{SW}_n\,2(1-n)\left( \frac{\partial}{\partial\vs} - \vs\calI\right)(\Hsp_{n}\expSW)
  =  \gamma^{SW}_n\,2(1-n)\left( \Hspp_{n} -\vs\Hsp_{n} - \vs\Hsp_{n}\right)\expSW\nonumber\\[0.5em]
  &= \gamma^{SW}_n\,2(1-n)\left( \Hspp_{n} - 2\vs\Hsp_{n} \right)\expSW
  =  \gamma^{SW}_n\,2(1-n)(-2n)\Hs_{n}\expSW = 4n(n-1)\psi_n .
\end{align}
The final case, for a generic $k$, follows by
a recursive argument, allowing us to prove that:
\begin{align}
  \Lst^k\Ls^k\psi_n = (-1)^k\,2^k\,\frac{n!}{(n-k)!}\psi_n, \qquad n\geq k.
\end{align}
Except for the factor $2^k$, this expression is the same as that
in~\eqref{eq:0076}.
Therefore, we conclude that every element of the SW Hermite functions'
basis is an eigenfunction of the $2k$-th operator $\Lst^{k}\Ls^{k}$
with eigenvalue $(-1)^{k}2^{k}\,n!\slash{(n-k)!}$, for $n\geq k$.  The
eigenvalue is zero for $0\leq n \leq k-1$.
We can similarly conclude that the $2k$-th operator does not modify
the equations for the first $k$ modes of the expansion of $\fs$.

We end this section by investigating the action of the
generalized Lenard-Bernstein operators on Hermite functions expressed as
linear combinations of SW Hermite basis functions.
To this purpose, we consider the expansion:
%$\fs=h\expSW$, where
%the polynomial function is given by
%% 
%\begin{align}
%  \hs(\vs) = \sum_{n=0}^\infty\Cs_{n}\Hs_{n}(\vs).
%  \label{eq:0070}
%\end{align}
%By multiplying and dividing by the normalization factor $\gamma^{AW}$,
%and, then, using the definition of the AW basis functions we find that
\begin{align}
  \fs 
  = \hs\expSW 
  = \left[\sum_{n=0}^\infty\Cs_{n}\Hs_{n}\right]\expSW
  = \sum_{n=0}^\infty\frac{\Cs_{n}}{\gamma^{SW}_n}\left[\gamma^{SW}_n\Hs_{n}\expSW\right]
  = \sum_{n=0}^\infty\Csb_{n}\psi_{n},
  \label{eq:0070b}
\end{align}
where $\Csb_{n}=\Cs_{n}\slash{\gamma^{SW}_{n}}$.
Since $\psi_{n}$ is an eigenfunction of the generalized Lenard-Bernstein
operators, we readily find the following relations:
\begin{align}
  \Lst\Ls\fs 
  &= \sum_{n=0}^{\infty}\Csb_{n}\Lst\Ls\psi_{n}
  = \sum_{n=0}^{\infty}(-2n)\Csb_{n}\psi_{n},\\[0.5em]
  %% --------------------------------------------
  \Lst^2\Ls^2\fs 
  &= \sum_{n=0}^{\infty}\Csb_{n}\Lst^2\Ls^2\psi_{n}
  = \sum_{n=0}^{\infty}4n(n-1)\Csb_{n}\psi_{n},\\[0.5em]
  %% --------------------------------------------
  \ldots\nonumber\\[0.5em]
  %% --------------------------------------------
  \Lst^k\Ls^k\fs 
  &= \sum_{n=0}^{\infty}\Csb_{n}\Lst^k\Ls^k\psi_{n}
  = \sum_{n=0}^{\infty} (-1)^{k}2^k\frac{n!}{(n-k)!}\Csb_{n}\psi_{n},
\end{align}
from which we deduce that:
\begin{align}
  \Lst\Ls\fs &= \sum_{n=0}^\infty\Ds^{(1)}_{n}\psi_{n}
  \qquad\textrm{with~}\Ds^{(1)}_{n} = -2n\Csb_{n},\\[0.5em]
  %% -----------------------------------------------------------
  \Lst^2\Ls^2\fs &= \sum_{n=0}^\infty\Ds^{(2)}_{n}\psi_{n}
  \qquad\textrm{with~}\Ds^{(2)}_{n} = 4n(n-1)\Csb_{n},\\[0.5em]
  %% -----------------------------------------------------------
  \ldots\nonumber\\[0.5em]
  %% -----------------------------------------------------------
  \Lst^{k}\Ls^{k}\fs &= \sum_{n=0}^\infty\Ds^{(k)}_{n}\psi_{n}
  \qquad\textrm{with~}\Ds^{(k)}_{n} = (-1)^{k}2^k\frac{n!}{(n-k)!}\Csb_{n}.
\end{align}
By definition, it holds that
$\Ds^{(k)}_{0}=\Ds^{(k)}_{1}=\ldots=\Ds^{(k)}_{k-1}=0$ for a generic
$k\geq1$.
%%
%%
% Consistently with the previous comments, we can set $\Ds^{(1)}_{0}=0$,
% for $k=1$; $\Ds^{(2)}_{0}=\Ds^{(2)}_{1}=0$, for $k=2$; and, for a
% generic $k>2$, we set
% $\Ds^{(k)}_{0}=\Ds^{(k)}_{1}=\ldots=\Ds^{(k)}_{k-1}=0$.
%%

\medskip 
As far as mass and momentum conservations are concerned, we do not
have the same results of the AW Hermite discretization.
Indeed, we can check that equations \eqref{eq:0125} and
\eqref{eq:0151} do not hold anymore in the symmetric case.
Instead, we can prove the conservation of the weighted integrals:
\begin{align*}
  \int_{\REAL}\fs(\vs,\ts)\expSW\ddv
  \quad\textrm{and}\quad
  \int_{\REAL}\vs\fs(\vs,\ts)\expSW\ddv,
\end{align*}
which however are not associated with physical, conserved quantities
of interest in the continuous setting.
%%
% However, differently from the AW Hermite case, the preservation of
% lower velocity moments in the SW Hermite case does not guarantee the
% existence of physically meaningful invariants such as the energy in
% the Vlasov model.

%%%%%%%%%%%%%%%%%%%%%%%%%%%%%%%%%%%%%%%%%%%%%%%%%%%%%%%%%%%%%%%%%%%%%%%%%%%%%%%

\smallskip
\section{Hermite approximations of the advection equation}
\label{sec:Hermite:advection}

We take into account the following time-dependent problem for the
unknown scalar field $\fs(\vs,t)$:
\begin{align}
  \frac{\partial\fs}{\partial t} - \frac{\partial\fs}{\partial\vs} = 0,
  \label{eq:0235a}
\end{align}
supplemented with the initial condition:
\begin{align}
  \fs(\vs,0) = \fs_0(\vs).
  \label{eq:0235}
\end{align}

% We will study here the stability of the Galerkin variational
% formulation either in the SW case or the AW case.
% %%
% %% \subsection{SW Hermite: Stability of $\fs$ in the weighted $\LTWO$ norm}
% %%
% To study the stability of the SW Hermite variational formulation of
% equation~\eqref{eq:0235a}, we set $\fs(\vs,t)=\hs(\vs,t)\expSW$ (where
% $\hs$ is a polynomial in $\vs$). 
%%
We start with the study of the stability of the SW Hermite variational
formulation of equation~\eqref{eq:0235a}.
To this end, we set $\fs(\vs,t)=\hs(\vs,t)\expSW$ (where $\hs$ is a
polynomial in $\vs$).
%%
%% \subsection{SW Hermite: Stability of $\fs$ in the weighted $\LTWO$ norm}
%%
Take $\fs$ as the test function, and integrate over $\REAL$.
We obtain:
\begin{align}
  0 
  = \int_{\REAL}\left(\frac{\partial\fs}{\partial\ts}-\frac{\partial\fs}{\partial\vs}\right)\fs\ddv
  = \int_{\REAL}\left( \frac{\partial}{\partial\ts}\left(\frac{\fs^2}{2}\right)-\frac{\partial}{\partial\vs}\left(\frac{\fs^2}{2}\right) \right)\ddv
  = \frac{1}{2}\frac{d}{\dt}\int_{\REAL}\hs^2\expAW\ddv,
\end{align}
since the integral of $\partial\fs^2\slash{\partial\vs}$ is zero
because $\fs(\vs,t)\to 0$ for $\vs\to\pm\infty$.
The relation above shows that the weighted norm of the function $\fs$,
solving equation~\eqref{eq:0235a} in weak form, is conserved (i.e. it does
not change in time).

The same is not going to be true for the AW case. 
In fact, we may try to study the stability with the same approach
followed before.
%%
%% as we will see in the coming subsection.
%%
%% \subsection{AW Hermite: Stability of $\fs$ in the weighted $\LTWO$ norm}
%%
% We attempt to study the stability of the AW Hermite variational
% formulation of equation~\eqref{eq:0235a} through the same approach
% followed before.
%%
This time we set $\fs(\vs,t)=\hs(\vs,t)\expAW$ (where $\hs$ is a
polynomial in $\vs$).
We then take $\hs$ as test function and integrate over $\REAL$.
We obtain:
\begin{align}
  0 
  = \int_{\REAL}\left(\frac{\partial\fs}{\partial\ts}-\frac{\partial\fs}{\partial\vs}\right)\hs\ddv
  = \int_{\REAL}\hs\frac{\partial\hs}{\partial\ts}\expAW\ddv
  - \int_{\REAL}\hs\frac{\partial\fs}{\partial\vs}\ddv.
\end{align}
Successively, we integrate by parts the last term, substitute
$\fs=\hs\expAW$ and integrate by parts again.
All the boundary terms are zero since they involve a polynomial in
$\vs$ multiplied by a decaying exponential and are omitted.
This procedure yields:
\begin{align}  
  - \int_{\REAL}\hs\frac{\partial\fs}{\partial\vs}\ddv
  = \int_{\REAL}\fs\frac{\partial\hs}{\partial\vs}\ddv
  = \int_{\REAL}\frac{\partial\hs}{\partial\vs}\,\hs\expAW\ddv
  = \int_{\REAL}\frac{\partial}{\partial\vs}\left(\frac{\hs^2}{2}\right)\expAW\ddv
  = \int_{\REAL}\hs^2\,v\,\expAW\ddv.
\end{align}
Finally, we find that:
\begin{align}
  0 
  = \frac{d}{\dt}\int_{\REAL}\frac{\hs^2}{2}\expAW\ddv
  + \int_{\REAL}\vs\,\hs^2\expAW\ddv.
\end{align}
Since $\vs\in\REAL$ can assume positive or negative values, the sign
of the second integral is undetermined, and therefore, the AW Hermite
variational formulation is not absolutely stable in the weighted
$\LTWO (\REAL)$ norm. 
Note, however, that the weighted norm in the AW case does not have a
direct physical meaning as in the SW case.
In both the continuous case and its SW Hermite discretization, the
quantity $\int_{\REAL}\fs^2 dv$ is preserved.
This quantity is not preserved in the AW discretization.
In fact, we are in the situation in which neither the weighted
$\LTWO$-norm nor the unweighted one are preserved.

Now, we derive the recursive equation for the coefficients of the
Hermite expansion in both AW and SW cases.
In order to simplify the notation, in the expressions below, we set
$\gamma_n=\gamma_{n}^{SW}$ when we deal with the SW case or
$\gamma_n=\gamma_{n}^{AW}$ when we deal with the AW case (we recall
that these coefficients are defined in~\eqref{eq:gamma_SW}
and~\eqref{eq:gamma_AW}).
Also, we use the notation $\Csb_{n}=\Cs_{n}\slash{\gamma_n}$ to denote
the coefficients of the expansion in the Hermite functions $\psi_n$.
As usual, we have:
\begin{align}
  \fs(\vs,t) 
  = \sum_{n=0}^{\infty}\Cs_{n}(t)H_n(\vs)\expAW
  = \sum_{n=0}^{\infty}\Csb_{n}(t)\psi_n(\vs).
  \label{eq:0240}
\end{align}
%where $\psi_n$ can be either the $n$-th AW or the $n$-th SW Hermite
%basis function introduced in~\eqref{eq:0245}.
%% 
Accordingly, the initial condition is set through the relation:
\begin{align}
  \sum_{n=0}^{\infty}\Cs_{n,0}H_n(\vs)\expAW =
  \sum_{n=0}^{\infty}\Csb_{n,0}\psi_n(\vs) = \fs_0(\vs).
 \label{eq:0235b}
\end{align}

To derive the system of equations for the coefficients $\Csb_n$
related to the solution of~\eqref{eq:0235a}, we
multiply~\eqref{eq:0235a} by $\psi^{m}$ and integrate in $\vs$ over
$\REAL$.
All integrals can easily be computed using the orthogonality of the
Hermite functions' basis.
In view of expansion~\eqref{eq:0240}, we have that:
\begin{align}
  0
  &= \int_{\REAL}\left(\frac{\partial\fs}{\partial\ts}-\frac{\partial\fs}{\partial\vs}\right)\psi^m\ddv 
  = \sum_{n=0}^{\infty}{\bigdot\Cs}^\star_n(t)\int_{\REAL}\psi_n\psi^m\ddv 
  - \sum_{n=0}^{\infty}\Csb_{n}(t)\int_{\REAL}\frac{d\psi_n}{d\vs}\psi^{m}\ddv
  \nonumber\\[0.5em]
  &= {\bigdot\Cs}^\star_m(t) - \sum_{n=0}^{\infty}\Csb_{n}(t)\int_{\REAL}\frac{d\psi_n}{d\vs}\psi^{m}\ddv,
  \label{eq:0258:SW}
\end{align}
where the upper dot indicates the derivative with respect to $t$.
The equation for each coefficient $\Csb_{n}(t)$ can be recovered by
reformulating $d\psi_n\slash{d\vs}$ in terms of the basis functions
$\psi_n$ and using the orthogonality against $\psi^{m}$.
We discuss the two cases for the AW and SW Hermite approximation in the
following subsections.

\subsection{Symmetrically-weighted case}
\label{subsec:SW-Hermite:advection}
To ease the notation in the developments of this section, we continue
using the symbol $\gamma_{n}$ instead of $\gamma_{n}^{SW}$, which is
defined in~\eqref{eq:gamma_SW}.
For $n\geq 1$, using~\eqref{eq:0hermite:2}-\eqref{eq:0hermite:3}, we
compute $d\psi_{n}\slash{d\vs}$ as follows:
\begin{align}
  \frac{d\psi_n}{\ddv} 
  &= \frac{d}{\ddv}\left(\gamma_n\Hs_{n}\expSW\right)                                             %\nonumber\\[0.5em]
  = \gamma_n\left( \Hs_{n}^{\prime} - \vs\Hs_{n} \right)\expSW                                      %\nonumber\\[0.5em]
  = \gamma_n\left( \frac{1}{2}\Hs_n^{\prime} + \frac{1}{2}\Hs_n^{\prime} - \vs\Hs_{n} \right)\expSW \nonumber\\[0.5em]
  &= \gamma_n\left( n\Hs_{n-1} - \frac{1}{2}\Hs_{n+1} \right) \expSW
  = \frac{n\gamma_{n}}{\gamma_{n-1}}\psi_{n-1}                        
  - \frac{1}{2}\frac{\gamma_{n}}{\gamma_{n+1}}\psi_{n+1}.
\end{align}
Thus, equation~\eqref{eq:0258:SW} implies that:
\begin{align}
  {\bigdot\Cs}^\star_{n}(t) 
  = \frac{(n+1)\gamma_{n+1}}{\gamma_{n}} \Csb_{n+1}(t) 
  - \frac{1}{2}\frac{\gamma_{n-1}}{\gamma_{n}}\Csb_{n-1}(t),
  \label{eq:0265}
\end{align}
that we supplement with the initial condition $\Csb_{n}(0)=\Csb_{n,0}$.
Equivalently, one has for $n\geq 1$:
\begin{align}
  \bigdot{\Cs}_{n}(t) 
  = (n+1)\Cs_{n+1}(t) - \frac{1}{2}\Cs_{n-1}(t),
  \label{eq:0265a}
\end{align}
with the (obvious) initial condition $\Cs_{n}(0)=\Cs_{n,0}$.
%%
% Recalling that
% $\gamma_{n}=\gamma_{n}^{SW}=\big(\sqrt{\pi}2^nn!\big)^{-1\slash{2}}$, we can
% make explicit the coefficients above:
% %%
% \begin{align}
%   \frac{\gamma_n}{\gamma_{n+1}}=\left(\frac{2^n n!}{2^{n+1}(n+1)!}\right)^{-1/2}
%   =\sqrt{2(n+1)},\qquad  \frac{\gamma_n}{\gamma_{n+1}}\frac{\gamma_n}{\gamma_{n-1}}=
%   \frac{\sqrt{2(n+1)}}{\sqrt{2n}}=\sqrt{\frac{n+1}{n}}.
%   \label{eq:0266}
% \end{align}
% This finally brings to a system of differential equations in the
% coefficients $\Cs_{n}$, $n\geq 1$, namely:
% \begin{align}
%   \bigdot{\Cs}_{n}(t) 
%   = \sqrt{\frac{n+1}{n}}\left( n\Cs_{n+1}(t) - \frac{1}{2}\Cs_{n-1}(t) \right), \qquad n\geq 1.  
%   \label{eq:0266a}
% \end{align}
The case $n=0$ can be treated separately, by observing that:
\begin{align*}
  \frac{d\psi_0}{\ddv} 
  &= \frac{d}{\ddv}\left( \gamma_0\Hs_0\expSW \right) 
  = -\gamma_0\vs\expSW 
  = -\frac{\gamma_0}{2\gamma_1}\left( \gamma_1\,2\vs\expSW \right) 
  = -\frac{\gamma_0}{2\gamma_1}\left( \gamma_1 \Hs_1\expSW \right) 
  \\[0.5em]
  &= -\frac{\gamma_0}{2\gamma_1}\psi_1 
  = -\frac{1}{\sqrt{2}}\psi_1 
  \quad \Rightarrow \quad \int_{\REAL}\frac{d\psi_0}{d\vs}\psi^{0}\ddv=0,
\end{align*}
since $\Hs_0(\vs)=1$, $\Hs_1(\vs)=2\vs$, and
$\gamma_0\slash{\gamma_1}=\sqrt{2}$, so obtaining
from~\eqref{eq:0258:SW} that
\begin{align}
  \bigdot{\Cs}_{0}(t) = 0 \quad \Rightarrow \quad C_0(t)=C_{0,0} \quad \forall t.
  \label{eq:0267}
\end{align}
We proved above that the system associated with
equations~\eqref{eq:0265}-\eqref{eq:0267} is stable
in the $\LTWO$-weighted norm.

\subsection{Asymmetrically-weighted case}
\label{subsec:AW-Hermite:advection}
As in the previous section we ease the notation by writing the symbol
$\gamma_{n}$ instead of $\gamma_{n}^{AW}$, which is defined
in~\eqref{eq:gamma_AW}.
In this case, using~\eqref{eq:0hermite:2}, multiplying and dividing by
$\gamma_{n+1}$, and using the definition of $\psi_{n+1}$, we have:
\begin{align}
  \frac{d\psi_n}{\ddv} 
  = \frac{d}{\ddv}\left[\gamma_n\Hs_{n}\expAW\right]
  = \gamma_n\left( \Hs_n^{\prime} - 2\vs\Hs_n \right)\expAW
  = -\gamma_n\Hs_{n+1}\expAW
  = -\frac{\gamma_n}{\gamma_{n+1}}\psi_{n+1},
  \label{eq:0270}
\end{align}
which now provides the differential equation, for $n\geq 1$:
\begin{align}
  {\bigdot\Cs}_{n}^\star (t) &= -\frac{\gamma_{n-1}}{\gamma_n}\Csb_{n-1}(t),
  \label{eq:0275}
\end{align}
supplemented with the initial condition $\Csb_{n}(0)=\Csb_{n,0}$.
This is equivalent to:
\begin{align}
  \bigdot{\Cs}_{n}(t) = - \Cs_{n-1}(t)
\label{eq:0277}
\end{align}
For $n=0$ we have again~\eqref{eq:0267}.
Moreover we have the initial conditions $\Cs_{n}(0)=\Cs_{n,0}$; hence,
$\Cs_{0}(t)=\Cs_{0,0}$ for every $t\geq0$.

We now provide a solution to such a system of equations. For instance,
when $n=1$, we need to solve:
\begin{align}
  {\bigdot\Cs}_{1}(t) = -C_{0}(t) \quad \Rightarrow \quad C_1(t)=C_{1,0}-C_{0,0} t.
  \label{eq:Cs:AW-Hermite}
\end{align}
Clearly, this coefficient grows in magnitude with $t$.
By successive integrations, one can prove that the $n$-th coefficient
behaves as $t^n$.
In practice, it is possible to find numbers $\alpha_\ell^{(n)}$ in
such a way that:
\begin{align}
  C_{n}(t) = \gamma_n  C_{n}^\star (t)=\sum_{\ell=0}^{n}\alpha_{\ell}^{(n)}t^{\ell},
  \label{eq:0279}
\end{align}
which is clearly unbounded for $t$ tending to infinity. 
We already proved that the Galerkin approximation of the advection
problem in the AW case is not unconditionally stable in the $L^2(\REAL)$-weighted norm.
For a polynomial of degree at most $N$, such a norm with respect to
$t$ is given by the sum $\Big(\sum_{n=0}^N\big(\Csb_n(t)\big)^2\Big)^{1/2}$.
A way to stabilize the approximation scheme is to introduce some
numerical dissipation. 
We note, however, that this may not be the only option.
%%
%% The result in~\eqref{eq:0279} says that there are no chances to get
%% stability, with respect to $t$, in a milder form, since the divergence
%% behavior affects each coefficient. 
%%
%% This is why it is imperative to introduce some dissipation. 
%%
We will study this problem in the next section.

\subsection{Some additional considerations on the SW and AW Hermite approximations}
We consider two exact solutions of equation~\eqref{eq:0235a} that are well-suited for the treatment with Hermite functions
(in the SW case and the AW case, respectively) and see how their expansion coefficients look like, in particular with respect 
to the time variable t. 
It has to be remarked, however, that the truncated series of an exact solution does not
coincide, in general, with the discrete solution obtained by the Galerkin process. So, the purpose of the following computation is only to illustrate why the approximations based on the SW or the AW Hermite functions may behave rather differently.

First, we consider the exact solution of~~\eqref{eq:0235a} given by
$\fs(\vs,\ts)=e^{-\frac{(\vs+\ts)^2}{2}}$ and denote its coefficients
with respect to the SW Hermite functions by $\Cs_n^{SW,ex}$, where the
superscript ``ex'' stands for ``exact''.
At $t=0$, only one coefficient is nonzero, i.e.,
$\Cs^{SW,ex}_0(0)=\pi^{\frac{1}{4}}$.
For a generic $t>0$, the expansion coefficients of $\fs$ are, for $n\geq0$:
\begin{align}
  \Cs^{SW,ex}_n(t) 
  = \int_{\REAL}e^{-\frac{(\vs+\ts)^2}{2}}\psi_n(\vs)\ddv
  = \gamma_{n}^{SW}e^{-\frac{t^2}{4}}\int_{\REAL}e^{-\left(\vs+\frac{t}{2}\right)^2}\Hs_{n}(\vs)\ddv
  = \sqrt{\pi}2^{n}\gamma_{n}^{SW}\left(-\frac{t}{2}\right)^ne^{-\frac{t^2}{4}},
  \label{eq:advection:example:SW}
\end{align}
where, for the integration, we used the convolution
formula~\cite{Gradshteyn-Ryzhik:2007}:
%% {\RED{Credo che la formula sia del Gradshteyn and Ryzhik}}
\begin{align}
  \int_{\REAL}e^{-(\xs-\ys)^2}\Hs_{n}(\xs)\dx = \sqrt{\pi}2^{n}\ys^{n}.
  \label{eq:exact:integration:formula}
\end{align}
Formula~\eqref{eq:advection:example:SW} shows that all the expansion
coefficients $\Cs^{SW,ex}_n(t)$ converge to zero for $t\to\infty$ including
the one with $n=0$ (note that the coefficient provided by the Galerkin
approximation, namely $\Cs_{0}(0)$, is instead constant in time).

For the AW case we consider the exact solution
of~\eqref{eq:0235a} given by $\fs(\vs,\ts)=e^{-(\vs+\ts)^2}$ and we similarly
denote its expansion coefficients as $\Cs_n^{AW,ex}(t)$.
The expansion of $\fs$ on the AW Hermite basis functions still
contains only one coefficient at the initial time $t=0$, i.e.,
$\Cs^{AW,ex}_0(0)=\sqrt{\pi}$.
The new coefficients look as follows:
\begin{align}
  \Cs^{AW,ex}_n(t) 
  = \int_{\REAL}e^{-(\vs+\ts)^2}\psi_n(\vs)\ddv
  = \sqrt{\pi}2^{n}\gamma_{n}^{AW}\left(-t\right)^n,
  \label{eq:advection:example:AW}
\end{align}
where, for the integration, we used again
formula~\eqref{eq:exact:integration:formula}.
Formula~\eqref{eq:advection:example:AW} shows that the expansion
coefficients $\Cs^{AW,ex}_n(t)$ diverge to $\pm\infty$ when $t\to\infty$, the
sign depending on $n$ being even or odd.
Of course, in these circumstances a remedy can be easily found by
introducing a shift in the Hermite basis as mentioned in the
introduction.
The fact that the expansion needs to be centered and rescaled
properly has been known for a long time but complicates the analysis and so it will be considered in future work.

\begin{comment}
\smallskip
%% 
\begin{remark}
  For a fixed $N$, the Galerkin method allows for the construction of
  the coefficient $C_N(t)$ in terms of the lower order coefficients
  (see~\eqref{eq:0277}). 
  %% 
  Such a procedure is incompatible with the request that
  $C_{N+1}(t)=0, \forall t$.
  %% 
  In fact, going backwards in the recursion formula~\eqref{eq:0277} we
  would find that all the coefficients must be identically zero.
\end{remark}
\end{comment}

\section{The advection equation with the stabilization term in the AW case}
\label{sec:adve:stab}
We start our analysis by adding the second-order ($k=1$) operator
$\nu\Lst\Ls$ to the right-hand side of the advection equation:
\begin{align}
  \frac{\partial\fs}{\partial t} - \frac{\partial\fs}{\partial\vs} = \nu\Lst\Ls\fs,
  \label{eq:0280}
\end{align}
which we solve for $\fs(\vs,\ts)$.
We will prove that the new term acts like a stabilization term.
%%
%% We start by observing that this is unfortunately not true in the
%% weighted $\LTWO (\REAL )$ norm for the AW case.
%%
To this end, we set $\fs=\hs\expAW$, take $\hs$ as the test function,
(we assume that $\hs$ is a polynomial in $\vs$ at every time), and
integrate~\eqref{eq:0280} over $\REAL$.
We substitute the stabilization term $\nu\Lst\Ls\fs$ with the
expression given in~\eqref{eq:0040} (or~\eqref{eq:0067} with $k=1$) to
obtain:
\begin{align}
  \int_{\REAL}
  \left(
    \frac{\partial\fs}{\partial\ts} 
    - \frac{\partial\fs}{\partial\vs}
  \right)\hs\ddv =
  \nu\int_{\REAL}\big(\Lst\Ls\fs\big)\,\hs\ddv =
  -\frac{\nu}{2}\int_{\REAL}\big(\hsp\big)^2\expAW\ddv.
  \label{eq:0283}
\end{align}
We integrate by parts the second integral term
%% $\hs(\partial\fs\slash{\partial\vs})$
and apply the Young inequality (with constant $\sigma$) to obtain:
\begin{align}
  \frac{1}{2}\frac{d}{\dt}
  \int_{\REAL}\hs^2\expAW\ddv
  &=
  -\int_{\REAL}\hsp\hs\expAW\ddv 
  -\frac{\nu}{2}\int_{\REAL}\big(\hsp\big)^2\expAW\ddv
  \nonumber\\[0.5em]
  &\leq 
  \ABS{\int_{\REAL}\hsp\hs\expAW\ddv}
  -\frac{\nu}{2}\int_{\REAL}\big(\hsp\big)^2\expAW\ddv
  \nonumber\\[0.5em]
  &\leq
  \frac{1}{2\sigma}\int_{\REAL}\hs^2\expAW\ddv +
  \frac{1}{2}(\sigma-\nu)\int_{\REAL}\big(\hsp\big)^2\expAW\ddv,
  \label{eq:AW:stability}
\end{align}
where we used the fact that the boundary contributions from the
integration by parts are zero.
From the Poincar\`e inequality~\eqref{eq:2503} (take
$\varphi=\hs$) we have that
\begin{align}
  -\frac{1}{2}\int_{\REAL}\big(\hs^{\prime}\big)^2\expAW\ddv
  \leq 
  -\int_{\REAL}\hs^2\expAW\ddv
  +\sqrt{\pi}\Cs_{0}^2.
  \label{eq:AW:stability:02}
\end{align}
Using this inequality with $\nu>\sigma$, we find that
\begin{align}
  \frac{1}{2}\frac{d}{\dt}
  \int_{\REAL}\hs^2\expAW\ddv
  \leq
  \left(\frac{1}{2\sigma}-(\nu-\sigma)\right)
  \int_{\REAL}\hs^2\expAW\ddv
  + (\nu-\sigma)\sqrt{\pi}\Cs_{0}^2.
  \label{eq:AW:stability:00}
\end{align}
The coefficient $\big(1\slash{(2\sigma)}-(\nu-\sigma)\big)$ is
negative if $\nu>\sigma+1\slash{(2\sigma)}$.
For example, by taking $\sigma=1$ and $\nu>3\slash{2}$, we find:
\begin{align}
  \frac{1}{2}\frac{d}{\dt}
  \int_{\REAL}\hs^2\expAW\ddv
  \leq
  -\left(\nu-\frac{3}{2}\right)\int_{\REAL}\hs^2\expAW\ddv
  + (\nu-1)\sqrt{\pi}\Cs_0^2.
  \label{eq:AW:stability:05}
\end{align}
Now, we consider $\Cs_{0}(t)=\Cs_{0}(0)=\Cs_{0,0}$ and introduce the quantities:
\begin{align}
  \Ks = \frac{\nu-1}{\nu-\frac{3}{2}}\sqrt{\pi}\Cs_{0,0}^2
  \qquad\textrm{and}\qquad
  \Ys(t) = \int_{\REAL}\hs^2\expAW\ddv - \Ks,
  \label{eq:Yt:def}
\end{align}
so we can rewrite~\eqref{eq:AW:stability:05} as
\begin{align}
  \frac{1}{2}\frac{d}{\dt}\Ys(t)\leq -\left(\nu-\frac{3}{2}\right)\Ys(t)
\end{align}
since $\Ks$ is constant.
Note that for $t=0$ we have
\begin{align}
  \Ys(0) = \int_{\REAL}\hs_{0}^2\ddv-\Ks,
  \label{eq:Y0:def}
\end{align}
where $\hs_{0}=\hs(\vs,0)$, which is provided by the expansion of the
initial solution $\fs_{0}$.
Finally, an application of the Gronwall's inequality yields
\begin{align}
  \Ys(t)
  \leq \Ys(0)\exp{\left(-2\Big(\nu-\frac{3}{2}\Big)t\right)}
  \leq \Ys(0),
  \label{eq:AW:stability:bound}
\end{align}
since the argument of the exponential is negative.
Using the expression of $\Ys(t)$ and $\Ys(0)$, respectively given
in~\eqref{eq:Yt:def} and~\eqref{eq:Y0:def}, the condition
$\Ys(t)\leq\Ys(0)$ implies that
\begin{align}
  \int_{\REAL}\hs^2\expAW\ddv
  %\leq \Ks + \int_{\REAL}\hs_{0}^2\ddv-\Ks
  \leq \int_{\REAL}\hs_{0}^2\ddv
  = \int_{\REAL}\hs(\vs,0)^2\ddv,
  \label{eq:AW:stability:08}
\end{align}
which is the stability in the weighted $L^2$ norm.
Note that $\nu>\frac{3}{2}$ is a sufficient but not necessary 
conditions for stability.

\medskip
Concerning the case $k>1$, a proof of stability for $\nu$ sufficiently large, can be given following
the same steps of the case for $k=1$.
We just provide here a sketch of the main steps for the classical $L^2$-weighted norm.
Thanks to~\eqref{eq:0067}, formula~\eqref{eq:0283} can be rewritten as
\begin{align}
  \int_{\REAL}
  \left(
    \frac{\partial\fs}{\partial\ts} 
    - \frac{\partial\fs}{\partial\vs}
  \right)\hs\ddv =-(-1)^{k}
  \nu\int_{\REAL}\big( \Lst^{(k)}\Ls^{(k)}\fs \big)\,\hs\ddv =
  -\frac{\nu}{2^k}\int_{\REAL}\big(\hs^{(k)}\big)^2\expAW\ddv.
  \label{eq:0290:AW:stab}
\end{align}
As in~\eqref{eq:AW:stability} we use the Schwarz and Young inequality;
then, we estimate the right-hand side of~\eqref{eq:0290:AW:stab} by
using~\eqref{eq:2503a} with $p=1$ and $m=k$.
By using~\eqref{eq:AW:stability:02}, we arrive at 
\begin{align}
  \frac{1}{2}\frac{d}{\dt}
  \int_{\REAL}\hs^2\expAW\ddv
  \leq
  \Phi_1\int_{\REAL}\hs^2\expAW\ddv + \Phi_2
  \label{eq:AW:stability:10}
\end{align}
where 
\begin{align}
  \Phi_1 = \frac{1}{2\sigma}-\nu(k-1)!+\sigma
  \quad\textrm{and}\quad
  \Phi_2 
  = \big(\nu(k-1)!-\sigma\big)\sqrt{\pi}\Cs_{0}^2
  + \nu(k-1)!\sqrt{\pi}\sum_{\ell=1}^{k-1}2^{\ell}\frac{(\ell!)^2}{(\ell-1)!}\Cs_{\ell}^2.
  \label{eq:AW:stability:15}
\end{align}
Now, we redefine
%\FNOTE{Ho spostato qui questo pezzo che prima stava alla fine della sezione 9}
\begin{align}
  \Ks = \frac{ \big(\nu(k-1)!-\sigma\big)}{ \nu(k-1)!-\sigma-\frac{1}{2\sigma} }\sqrt{\pi}\Cs_{0}^2
  \qquad\textrm{and}\qquad  
  \Ys(t) = \int_{\REAL}\hs_{0}^2\ddv-\Ks,
\end{align}
so that 
\begin{align}
  \frac{1}{2}\frac{d}{\dt}\Ys(\ts)\leq \Phi_1\Ys(\ts) + \Psi_1(t),
  \qquad\textrm{where}\qquad
  \Psi_1(\ts) = \nu(k-1)!\sqrt{\pi}\sum_{\ell=1}^{k-1}2^{\ell}\frac{(\ell!)^2}{(\ell-1)!}\Cs_{\ell}^2,
\end{align}
since $\Cs_{0}=\Cs_{0,0}$ is independent of $\ts$.
An application of the Gronwall's Lemma leads to
\begin{align*}
  \Ys(t)
  %\leq \left[\Ys(0)+\int_{0}^{\ts}\Psi_1(\tau)e^{ 2\Phi_1\tau }\ds\tau\right]e^{ -2\Phi_1\ts }
  \leq \Ys(0)e^{ -2\Phi_1\ts } + \int_{0}^{\ts}\Psi_1(\tau)\ds\tau 
\end{align*}

Choosing, for example, $\sigma=1$ and taking $\nu(k-1)!>3\slash{2}$,
it is easy now to get the stability estimate that generalizes~\eqref{eq:AW:stability:bound}
to any $k\geq1$.
We also note that the diffusion parameter $\nu$ is now multiplied by $(k-1)!$.
So, if we increase $k$, the numerical diffusion due to the Lenard-Bernstein operators acts 
only on higher terms in the expansion of $\fs$ and we may probably take smaller values for
$\nu$.
%%
%% However, the first $k$ expansion coefficients from $0$ to $\ks$ are not controlled
%% by the above inequality.
%%
%textbf{[GM: vogliamo aggiungere altri commenti? Oppure modificare il
%commento originale:} \it If the last condition on $\nu$ is not
%satisfied, we expect that t is not possible to prove the stability in
%the classical $L^2$-weighted norm.  ote also that a suitable weighted
%norm as the one discussed before may e used to prove stability in a
%milder form.

\medskip
We confirm the stability result for $k=1$ by deriving the explicit recursive
formula for the Hermite expansion coefficients and providing their
explicit form.
To this end, we consider the second expansion of $\fs$ given
in~\eqref{eq:0240} and repeat the calculation of
Section~\ref{subsec:AW-Hermite:advection}
%% the previous sections
by including now the stabilization term $\nu\Lst\Ls\fs$, which can be
treated in the AW case with the help of~\eqref{eq:0030a}:
\begin{align}
  0
  &= 
  \int_{\REAL}
  \left(
      \frac{\partial\fs}{\partial\ts}
    - \frac{\partial\fs}{\partial\vs}
    - \nu\Lst\Ls\fs
  \right)\psi^m\ddv 
  \nonumber\\[0.5em]
  &
  = \sum_{n=0}^{\infty}\bigdot{C}_{n}^\star (t)\int_{\REAL}\psi_n\psi^m\ddv 
  - \sum_{n=0}^{\infty}\Csb_{n}(t)\int_{\REAL}\frac{d\psi_n}{d\vs}\psi^{m}\ddv
  + \nu\sum_{n=0}^{\infty}n\Csb_{n}(t)\int_{\REAL}\psi_n\psi^{m}\ddv
  \nonumber\\[0.5em]
  &= {\bigdot C}_{m}^\star (t) 
  - \sum_{n=0}^{\infty}\Csb_n(t)\int_{\REAL}\frac{d\psi_n}{d\vs}\psi^{m}\ddv
  + \nu m\Csb_m(t).
  \label{eq:0258}
\end{align}
We compute the last integral using again~\eqref{eq:0270} to obtain:
\begin{align}
  {\bigdot C}_{n}^\star(t) = -\frac{\gamma_{n-1}}{\gamma_{n}}\Csb_{n-1}(t) - \nu n\Csb_{n}(t),
  \label{eq:0290}
\end{align}
which holds for $n\geq 1$, while for $n=0$ we find that
$\Csb_{0}(t)=\Csb_{0,0}$ is constant.
%% 
% \begin{align}
%   \gamma_{n} := \gamma_{n}^{AW} = (\sqrt{\pi}2^n\,n!)^{-\frac{1}{2}}.
%   %%\widetilde{\gamma}_{n}^{AW} = (2^n\,n!)^{-\frac{1}{2}}.
% \end{align}
% 
We rewrite the above system of equations as follows (compare
with~\eqref{eq:0277}):
\begin{align}
  \bigdot{\Cs}_{n}(t) = -\Cs_{n-1}(t) - \nu n\Cs_{n}(t).
  \label{eq:0290a}
\end{align}

For instance, for $n=1$, we find the ordinary differential equation:
%\begin{align}
%  \bigdot{\Csb}_{1}(t) = 
%  -\frac{\gamma_1}{\gamma_2} \Csb_{0}(t) 
%  -\frac{2\nu}{\gamma_1}\Csb_{1}(t),
%  \label{eq:0295}
%\end{align}
%%
%or, moving back to the coefficient $\Cs_1$ through~\eqref{eq:0290a}:
\begin{align}
  \bigdot{\Cs}_{1}(t) = -\Cs_{0}(t) - \nu\Cs_{1}(t) = -C_{0,0} - \nu\Cs_{1}(t),
  \label{eq:0295}
\end{align}
% and
% \begin{align}
%   \Csb_{1}(0) = \Csb_{1,0} = \gamma_{1}\int_{\REAL}\fs_{0}(\vs)\Hs_{1}(\vs)\ddv.
% \end{align}
%%
the solution of which is:
\begin{align}
  C_{1}(t) = \left(C_{1,0}+\frac{1}{\nu}C_{0,0}\right)e^{-\nu t}-\frac{1}{\nu}C_{0,0}.
  \label{eq:0296}
\end{align}
Since $\nu$ is positive, $C_1(t)$ is clearly bounded with respect to
$t$.
%
%\begin{align}
%  \Csb_{1}(t) = \frac{1}{\beta}\left( e^{\beta t+\epsilon} - \alpha \right),
%\end{align}
%where
%\begin{align}
%  \alpha = -\frac{\gamma_1}{\gamma_2}\Csb_{0,0} = -2\Csb_{0,0},
%  \quad
%  \beta = -\frac{2\nu}{\gamma_1} = -2\sqrt{2\pi}\nu,
%\end{align}
%and $\epsilon$ is such that 
%\begin{align}
%  \Csb_{1}(0) = \Csb_{1,0} := \int_{\REAL}\fs(\vs,0)\psi^{n}(\vs)\ddv.
%\end{align}
%% 
%In terms of $\Cs_1(t)$:
%\begin{align}
%  \Cs_1(t) = -\big(e^{-\nu t+c}+\sqrt{2}\big)\slash{\nu}
%\end{align}
%and $c = \textrm{ln}( -\nu\Cs_{1,0}-\sqrt{2} )$, which follows by
%imposing that $\Cs_{1}(0)=\Cs_{1,0}$.

It is not hard to show that, for a generic $n$, the expression of the
$n$-th coefficient takes the form:
\begin{align}
 C_n(t) = \sum_{\ell=0}^{n}\alpha_{\ell}^{(n)}e^{-\ell \nu t},
 \label{eq:0299}
\end{align}
where the constants $\alpha_{\ell}^{(n)}$ depend on $n$ and the
diffusion parameter $\nu$.
It is important to analyze such a dependence on the
diffusion parameter.
Indeed, using~\eqref{eq:0299} in~\eqref{eq:0290a} for $n\geq1$ and
  $0\leq\ell\leq n-1$ yields the recursive relation
\begin{align*}
  \alpha_{\ell}^{(n)} = -\frac{1}{\nu(n-\ell)}\alpha_{\ell}^{(n-1)},
\end{align*}
from which a straightforward calculation yields:
\begin{align*}
  \alpha_{\ell}^{(n)} = \frac{(-1)^{n-\ell}}{(n-\ell)!\,\nu^{n-\ell}}\alpha_{\ell}^{(\ell)}.
\end{align*}
%% where for $\alpha_{0}^{0}=\Cs_{0,0}$.
%%
From the initial condition
\begin{align*}
  \Cs_{n,0}
  = \Cs_{n}(0)
  = \sum_{\ell=0}^{n}\alpha_{\ell}^{(n)}
  = \alpha_{n}^{(n)} + \sum_{\ell=0}^{n-1}\alpha_{\ell}^{(n)},
\end{align*}
we find the expression of $\alpha_{n}^{(n)}$, which is given by
\begin{align*}
  \alpha_{n}^{(n)} 
  = \Cs_{n,0} - \sum_{\ell=0}^{n-1}\alpha_{\ell}^{(n)}
  = \Cs_{n,0} 
  + \sum_{\ell=0}^{n-1}\frac{(-1)^{n-\ell}}{ (n-\ell)!\,\nu^{n-\ell} }\alpha_{\ell}^{(\ell)}.
\end{align*}
For example, starting from $\alpha_{0}^{0}=\Cs_{0,0}$, for $n=1$, we
find that
$\alpha_{1}^{(1)}=\Cs_{1,0}-\alpha_{0}^{0}\slash{\nu}=\Cs_{1,0}-\Cs_{0,0}\slash{\nu}$.
Similarly, $\alpha_{2}^{(2)}$ is computed from $\alpha_{0}^{(0)}$ and
$\alpha_{1}^{(1)}$, and the following coefficients are
obtained from the ones that have already been computed.  
One can realize that $\nu$ appears at the
denominator to the $n$-th power.
%% (see for instance~\eqref{eq:0296}, corresponding to $n=1$).
%%
It turns out that the coefficients
$\Cs_{n}(t)$ in~\eqref{eq:0299} are of the form $\Cs_n(0)$ plus a
dissipative term.
The stronger dissipation is obtained when $\ell=1$, which provides a
contribution like $e^{-\nu t}\slash{\nu}$ (see~\eqref{eq:0296}).
If we do not want this dissipation to be too heavy so that the
perturbation is of order $\varepsilon$ when we integrate until the
final time $T$, we can consider 
$e^{-\nu\Ts}\approx\nu\varepsilon$ and take 
$T\approx\ABS{\ln(\nu\epsilon)}\slash{\nu}$.

\section{Time discretization of the 1-D problem}
\label{sec:time:discretization}
We study the numerical approximation of the
system of differential equations in~\eqref{eq:0290a}.
We use an implicit conservative method in time such as the trapezoidal
rule.
For a time-step $\Delta t>0$, we write for $j\geq 1$:
\begin{align}
  \frac{\Cs_n^j - \Cs_n^{j-1}}{\Delta t} = -\frac{\Cs_{n-1}^j + \Cs_{n-1}^{j-1}}{2}  
  - \nu n\frac{ \Cs_n^j + \Cs_n^{j-1}}{2},
  \label{eq:0290d}
\end{align}
with the initial condition $\Cs_n^0=C_{n,0}$. 
For $n=0$ we have instead $\Cs_0^j=\Cs_{0,0}, \forall j\geq 0$.
For instance, we can make the formula explicit for $n=1$:
\begin{align}
  \Cs^j_1\left( 1+\frac{\nu}{2}\Delta t\right) 
  = \Cs^{j-1}_1\left( 1-\frac{\nu}{2}\Delta t\right) -\Delta t\Cs_{0,0}.
\end{align}
After defining $\chi_n=(1-\frac12 \nu n \Delta t)/(1+\frac12 \nu n
\Delta t), n\geq 1$, we get $\vert \chi_n\vert <1$.
By recursive arguments, one can show that the expression for $\Cs_1^j$
takes the form of a linear combination of powers of $\chi_1$, i.e.:
\begin{align}
  \Cs^j_1= \sum_{\ell=0}^{j} (\chi_1)^\ell \alpha_\ell,
\end{align}
where the numbers $\alpha_\ell$ depend on $\nu$ and $\Delta t$.
This expression is inserted in~\eqref{eq:0290d} in order to compute
the sequence $\Cs_2^j, \forall j\geq 0$, and so on.

We may assume that the solution $h=fe^{v^2}$ of~\eqref{eq:0280}
belongs to the space of polynomials of degree less or equal to $N$.
When $n$ reaches the value $N$, the expression of the corresponding
coefficients $\Cs_N^j, \forall j\geq 0$ is a combination of all the
powers $(\chi_n)^\ell$ with $1\leq n\leq N$ and $0\leq \ell\leq j$.

Since $\vert\chi_N\vert <1$, the discretization method is always
unconditionally stable.
However, a wise relation between the parameters $N$, $\nu$ and $\Delta
t$ should be set up in order to avoid unpleasant numerical effects due
to the {\sl stiffness} of the originating differential
system~\eqref{eq:0290a} for $N$ large.
A rule of thumb is to require that the product $\nu N \Delta t$ is of the order
of unity.
Actually, if we analyze~\eqref{eq:0299} when $n=N$, the most
significant term is that given by the exponential $e^{-N\nu t}$,
displaying a very steep tangent for $t=0$.
Although there are in principle no restrictions on $\Delta t$ for the
trapezoidal scheme, such quick variations in time are well resolved only
if the time-step is maintained suitably small.

\medskip

The last arguments show that stability holds for any $\nu >0$, whereas 
in~\eqref{eq:AW:stability:08} the proof was only provided for $\nu>3/2$.
Indeed, we conjecture that the stability in the $L^2$-weighted norm is
not verified for values of $\nu$ less than a certain constant.
However, it is possible to construct milder weighted norms where a
result of stability can still be achieved for any $\nu$.
We show how to do this by starting from the differential
system~\eqref{eq:0290a}.
For any $n\geq 1$, we multiply~\eqref{eq:0290a} by the Hermite coefficient $\Cs_n$
  and use the Young inequality on the right-hand side to obtain
\begin{align}
  \frac{1}{2}\frac{d}{\dt}\Cs_{n}^2 
  = -\Cs_{n}\Cs_{n-1} - \nu n\Cs_{n}^2
  \leq \frac{1}{2\sigma_n}\Cs_{n}^2 + \frac{\sigma_n}{2}\Cs_{n-1}^2 - \nu n\Cs_{n}^2.
\end{align}
The family of parameters $\sigma_{n}>0$ will be decided later on.
We multiply both sides of the inequality above by a weight $\ws_{n}>0$
and sum over index $n$, so obtaining
\begin{align}
  \frac{1}{2}\frac{d}{\dt}\sum_{n=1}^{\infty}\ws_{n}\Cs_{n}^2
  \leq 
  \sum_{n=1}^{\infty}\frac{1}{2\sigma_n}\ws_{n}\Cs_{n}^2 +
  \sum_{n=1}^{\infty}\frac{\sigma_{n}}{2}\ws_{n}\Cs_{n-1}^2 -
  \sum_{n=1}^{\infty}\nu n\ws_{n}\Cs_{n}^2.
\end{align}
By shifting the index in the sum containing $\Cs_{n-1}$ and collecting
the corresponding terms under the same symbol of summation, we get
\begin{align}
  \frac{1}{2}\frac{d}{\dt}\sum_{n=1}^{\infty}\ws_{n}\Cs_{n}^2
  \leq\sum_{n=1}^{\infty}\bigg[\Big(\frac{1}{2\sigma_{n}}-\nu n\Big)\ws_{n}+\frac{\sigma_{n+1}}{2}\ws_{n+1}\bigg]\Cs_{n}^2
  + \frac{\sigma_1}{2}\ws_{1}\Cs_{0}^2.
\end{align}
For example, we can consider to choose $\sigma_{n}=1\slash{(\nu n)}$.
Successively, we impose that the expression in the square brackets is
equal to $-(\nu /4)\ws_{n}$, which implies that
\begin{align}
  \ws_{n+1} = \nu^2(n+1)\left(n-\frac{1}{2}\right)\ws_{n},
  \label{eq:wn:recursive:def}
\end{align}
and starting for example from $\ws_{1}=1$ provides the full weight
sequence (altough a different starting value produce a different
sequence of coefficients $\ws_{n}$, our argument is independent of
such starting):

\begin{align}
  \frac{1}{2}\frac{d}{\dt}\sum_{n=1}^{\infty}\ws_{n}\Cs_{n}^2
  \leq -\frac{\nu}{4}\sum_{n=1}^{\infty}\ws_{n}\Cs_{n}^2 + \frac{1}{2\nu}\ws_{1}\Cs_{0}^2.
  \label{eq:dYdt:00}
\end{align}
Finally, by setting
\begin{align}
  \Ys(t) 
  = \sum_{n=1}^{\infty}\ws_{n}\Cs_{n}^2 - \big(2\slash{\nu}^2\big)w_1\Cs_{0}^2,
  \label{eq:Y:def}
\end{align}
we obtain 
\begin{align}
  \frac{1}{2}\Ys'(t) \leq -\frac{\nu}{4}\Ys(t).
  \label{eq:dYdt:10}
\end{align}
Thus, by applying the Gronwall's lemma, we conclude with the estimate
\begin{align}
  \Ys(t) \leq \Ys(0) e^{-(\nu /2)t}\leq \Ys(0),
\end{align}
for all $t\geq 0$,   
%%
%\begin{align*}
%  -\frac{\nu}{4}\sum_{n=1}^{\infty}\ws_{n}\Cs_{n}^2 + \frac{1}{2\nu}\ws_{1}\Cs_{0}^2
%  \leq -\frac{\nu}{4}\Ys(t)
%  =  -\frac{\nu}{4}
%  \left[ 
%    \sum_{n=1}^{\infty}\ws_{n}\Cs_{n}^2 + \big(2\slash{\nu}^2\big)w_1\Cs_{0}^2
%  \right].
%\end{align*}
%\textbf{
%Se fosse vero, troveremmo che}
%\begin{align*}
%  \frac{1}{2\nu}\ws_{1}\Cs_{0}^2 
%  \leq -\frac{\nu}{4}\left( \big(2\slash{\nu}^2\big)w_1\Cs_{0}^2 \right)
%  = -\frac{1}{2\nu}\ws_{1}\Cs_{0}^2
%\end{align*}
%\textbf{
%che \`e ovviamente falso (escluso il caso $0\leq0$, ovviamente).
%%% 
%Proporrei di riscrivere il pezzo sopra come segue.}
%\medskip
%\noindent
%% 
%In this way, we arrive at the inequality
%\begin{align*}
%  \frac{1}{2}\frac{d}{\dt}\sum_{n=1}^{\infty}\ws_{n}\Cs_{n}^2
%  \leq -\frac{\nu}{4}\sum_{n=1}^{\infty}\ws_{n}\Cs_{n}^2 + \frac{1}{2\nu}\ws_{1}\Cs_{0}^2.
%  %\label{eq:dYdt:00}
%\end{align*}
%By setting
%\begin{align*}
%  \widetilde{\Ys}(t) 
%  = \sum_{n=1}^{\infty}\ws_{n}\Cs_{n}^2 - \big(2\slash{\nu}^2\big)w_1\Cs_{0}^2,
%  %\label{eq:Y:def}
%\end{align*}
%we find out that 
%\begin{align*}
%  \frac{1}{2}\widetilde{\Ys}'(t) \leq -\frac{\nu}{4}\widetilde{\Ys}(t).
%  %\label{eq:dYdt:10}
%\end{align*}
%Thus, by applying the Gronwall's lemma we conclude with the estimate
%\begin{align*}
%  \widetilde{\Ys}(t) \leq \widetilde{\Ys}(0) e^{-(\nu /2)t}\leq \widetilde{\Ys}(0),
%\end{align*}
from which we can find our stability result 
\begin{align*}
  \sum_{n=1}^{\infty}\ws_{n}\Cs_{n}^2(t) \leq \sum_{n=1}^{\infty}\ws_{n}\Cs_{n}^2(0),
\end{align*}
since the term $\big(2\slash{\nu}^2\big)w_1\Cs_{0}^2$ in~\eqref{eq:Y:def} is independent of time and can be removed.
%%
%Finally, by adding the constant term above and setting
%\begin{align*}
%  \Ys(t) 
%  = \sum_{n=1}^{\infty}\ws_{n}\Cs_{n}^2 + \big(2\slash{\nu}^2\big)w_1\Cs_{0}^2,
%  %\label{eq:Y:def}
%\end{align*}
%we find that  
%\begin{align*}
%  \Ys(t)\leq\Ys(0),
%\end{align*}
%for all $t\geq 0$, which is our stability result.

The next step is to characterize the weights $w_n$ which are required
to satisfy the recursive relation~\eqref{eq:wn:recursive:def}.
Assuming that $\ws_{1}=1$, from a straightforward calculation, we find 
\begin{align}
  \ws_{n} 
  &= \big(2\nu^2\big)^{n-1}\,n!\,\frac{(2n-3)!}{2^{n-2}\,(n-2)!}
  = 2\big(\nu^2\big)^{n-1}\,n(n-1)\,(2n-3)!
\end{align}
By substituting into~\eqref{eq:Y:def}, we are finally able to give an expression to the
stability norm. 
Note that it depends on $\nu$.
We can go ahead with our computations by noting that
\begin{align}
  w_n &\geq 2\big(\nu^2\big)^{n-1}\,n(n-1)2^{n-2}\,(n-2)!
  =\frac12\big(\nu^2\big)^{n-1}\,2^{n}\,n!
\end{align}
%\begin{align}
%  \norm{\fs}{L^2_w} = \sqrt{\pi}\sum_{n=0}^{\infty}2^n \,n!\,\Cs_{n}^2.
%\end{align}
\noindent
Therefore, if $\nu\geq 1$, and, hence, $\nu^{2n}\geq1$ , we
discover that
\begin{align}
  \Ys(t) 
  &= \sum_{n=1}^{\infty}\ws_{n}\Cs_{n}^2 + \frac{2}{\nu^2}\Cs_{0}^2
  \geq \frac{1}{2}\sum_{n=1}^{\infty}(\nu^2)^{n-1}2^n\,n!\,\Cs_{n}^2 + \frac{2}{\nu^2}\Cs_{0}^2
  \geq \frac{1}{2\nu^2}\sum_{n=1}^{\infty} \nu^{2n} 2^{n}\,n!\Cs_{n}^2 + \frac{1}{2\nu^2}\Cs_{0}^2
  \nonumber\\[0.5em]
  &{\geq} \frac{1}{2\sqrt{\pi}\nu^2}\left(\sqrt{\pi}\sum_{n=0}^{\infty}2^n \,n!\,\Cs_{n}^2\right).
%	\norm{\fs}{L^2_w}.
\end{align}
Thus, if we can bound $Y$, we automatically bound the last term in
parenthesis, which corresponds to the square of the classical
$L^2$-weighted norm of the solution $\fs$ expanded as
in~\eqref{eq:0240}.
This confirms that, if $\nu$ is sufficiently large, stability is
ensured in the standard way.
On the other hand, when $\nu<1$, we can only rely on the stability
result involving the weights $w_n$.

\medskip
\noindent
If we are in finite dimension ($n\leq N$), the norms are equivalent
for any $\nu >0$, but with constants heavily dependent on $N$.
For example, for $\nu\leq 1$, which implies that
$(\nu^2)^{n-1}\geq(\nu^2)^{N-1}$, we can write
\begin{align}
  \Ys(t) = 
  \sum_{n=1}^{N}\ws_{n}\Cs_{n}^2 + \frac{2}{\nu^2}\Cs_{0}^2
  %%\geq \frac{1}{2}\sum_{n=1}^{N}\big(\nu^2\big)^{n-1}\,2^{n}\,n!\Cs_{n}^2 + \frac{2}{\nu^2}\Cs_{0}^2
  \geq\frac{\big(\nu^2\big)^{N-1}}{2}\sum_{n=1}^{N}2^n\,n!\,\Cs_{n}^{2} + \frac{2}{\nu^2}\Cs_{0}^2
  %% \nonumber\\[0.5em]
  \geq \frac{\big(\nu^2\big)^{N-1}}{2\sqrt{\pi}}\left(\sqrt{\pi}\sum_{n=0}^{N}2^n \,n!\,\Cs_{n}^2\right).
  \label{eq:stability:control}
\end{align}
This shows that, when $Y$ is bounded by a constant, the classical
$L^2$-weighted norm of the solution $\fs$ 
is bounded by that constant multiplied by a factor
behaving as the inverse of $\nu^{2N}$.
If we choose $\nu<1$, such a constant grows to infinity as $\mathcal{O}(\nu^{2N})$, 
and the stability control on the $L^2$-weighted norm of the solution $\fs$ provided
by inequality~\eqref{eq:stability:control} is lost.

% ----------------------------------------------------------------------------------

\section{Full discretization of the Vlasov-Poisson equation}
\label{sec:full-discretization:Vlasov-Poisson}
\setcounter{numbs}{6}

%%Let $\Omv=\REAL$, 
%%
We consider the AW Hermite-based discretization of the Vlasov-Poisson
problem~\eqref{eq:Vlasov:1}-\eqref{eq:Vlasov:2} for the
distribution function $\fs(\xs,\vs,t)=\hs(\xs,\vs,t)\expAW$ stabilized
by the Lenard-Bernstein-like operator of order $2k$ with $k\geq1$,
which we rewrite here for convenience of exposition:
\begin{align}
  &\frac{\partial\fs}{\partial t} + \vs\frac{\partial\fs}{\partial x} 
  -\Es\frac{\partial\fs}{\partial\vs} = 
  -(-1)^k \nu\Lst^{(k)}\Ls^{(k)}\fs 
  \qquad\textrm{in~}\Omega\times[0,T],
  \label{eq:1205:a}
  \\[0.5em]
  &\frac{\partial\Es}{\partial x} = 1 - \int_{\Omv}\fs\ddv 
  \qquad\textrm{in~}\Omega\times[0,T].
  \label{eq:1205:b}
\end{align}

System~\eqref{eq:1205:a}-\eqref{eq:1205:b} is completed by assigning
a sufficiently regular initial solution
$\fs(\xv,\vv,0)=\fs_0(\xs,\vs)$. 
We specialize the discussion to periodic boundary conditions in space,
i.e., at the boundaries of $\Omega_x$.

\medskip
Some of the reasons for approaching the Vlasov problem by Hermite discretizations 
have been pointed out in the introduction. 
The AW context is the one that guarantees a large number of conservation 
properties, even with the addition of the diffusion term discussed so far. 
By the way, from the practical viewpoint the use of the viscous term 
should not just be interpreted as a way to improve the time-stability 
of the schemes, but has an important role in the reduction of the negative 
phenomenon known as \textit{filamentation},
cf.~\cite{Camporeale-Delzanno-Bergen-Moulton:2015}, which shows up as a polluting 
effect on the computed solutions, due to the nonlinearity of the problem 
in conjunction with the truncation of the high modes.

%\RED{
%Concerning equation~\eqref{eq:0235a} of the next section, unweighted
%mass conservation has been proved for the SW Hermite case
%\cite{Holloway:1996,Schumer-Holloway:1998}.
%%
%Nevertheless, according to previous considerations, the successive
%advection-diffusion equation~\eqref{eq:0280} for the SW case does
%conserve neither the weighted nor the unweighted quantities.
%%
%On the other hand, as it will be made clear in the subsequent
%sections, the SW Hermite discretization of equation~\eqref{eq:0235a}
%is stable without the diffusive term.
%}

\medskip
To discretize the Vlasov-Poisson equations in time, we
integrate equation~\eqref{eq:1205:a} 
with respect to the independent unknown $\ts$ 
between $\ts^{j-1}$ and $\ts^{j}$ by applying 
the trapezoidal rule and we evaluate equation~\eqref{eq:1205:b} 
at $\ts^{j}$.
To ease the exposition, we assume a constant time step
$\Delta\ts=\ts^{j}-\ts^{j-1}$.
At the timestep $j\geq 1$, system~\eqref{eq:1205:a}-\eqref{eq:1205:b}
yields
\begin{align}
  \frac{\fs^{j}-\fs^{j-1}}{\Dt} 
  &+ \vs\frac{\partial}{\partial x}\left( \frac{\fs^{j}+\fs^{j-1}}{2} \right)
  - \frac{\Es^{j}+\Es^{j-1}}{2}\,\frac{\partial}{\partial\vs}
  \left( \frac{\fs^{j}+\fs^{j-1}}{2} \right)\nonumber\\[0.5em]
  &= -(-1)^k\nu\Lst^{(k)}\Ls^{(k)}\left( \frac{\fs^{j}+\fs^{j-1}}{2} \right) 
  \label{eq:1215}
  \\[0.5em]
%\end{align}
%\begin{align}
  \frac{\partial \Es^{j}}{\partial x} &= 1 - \int_{\Omv}\fs^{j}\ddv.
  \label{eq:1216}
\end{align}
For $j=0$ we impose the value of $f$ at time $t=0$ as initial datum.

Following the guidelines of the previous section, a proof of the absolute
stability in time of this scheme can be provided for a sufficiently large
parameter $\nu$.
The situation gets more technically involved if $\nu$ is relatively
small.
We remind you that in Section~\ref{sec:adve:stab} we distinguished
between $\nu\geq1$ and $\nu<1$.
In the latter case, stability is achieved in a suitable norm and the
generalization of this proof to the Vlasov-Poisson system becomes rather
complicated.

Here, our goal is to derive stability conditions that relate the time step
$\Delta t$, the collisional factor $\nu$ and the degree of the Hermite 
polynomial $N$.
%%
% regardless of a complete convergence analysis, which is out of the
% scope of this work.
%%
% We know from the previous sections that there are little chances to
% get a stability result in the $L^2(\Omega )$ weighted norm in the AW
% case.
%%
%% By the way, we can try a preliminary analysis on the well-posedness
%% of~\eqref{eq:1215} in favour of the unknown $\fs^{j}$.
%%
To this end, we write~\eqref{eq:1215} in operator form by collecting
all the terms involving the unknown variable $\fs^{j}$ on the
left-hand side and denoting all other terms that are computable from
what is known from the previous time step in the right-hand side term
$\gs^{j-1}$:
\begin{align}
  &\left[ \calI 
    + \frac{\Dt}{2}\vs\frac{\partial}{\partial x}
    - \frac{\Dt}{2}\left(\frac{\Es^{j}+\Es^{j-1}}{2}\right)\frac{\partial}{\partial\vs} 
    +(-1)^k \frac{\Dt}{2}\nu\Lst^{(k)}\Ls^{(k)}
  \right] \fs^{j} = \gs^{j-1}.
  \label{eq:1220}
  % \\[0.5em]
  % &\int_{\Omv}\fs^{k+1}\ddv + \frac{\partial\Es^{k+1}}{\partial x} = 1,
  % \label{eq:1225}
\end{align}
In this preliminary analysis, we will not take into consideration that
the problem is actually nonlinear.
Indeed, the value $\Es^j$ has still to be computed, since it is strictly
linked to $\fs^j$ through the relation~\eqref{eq:1216}. 
%%
%% We omit however the discussion of this aspect.

We first set $\fs^j=\hs^j\expAW$.
To simplify the exposition, we remove the label $j$ from $\hs^j$
and introduce the notation:
\begin{align}
  &\calA(\xs)  = \left(\int_{\Omv}\hs^2\expAW\ddv\right)^{\frac{1}{2}},\qquad
  \calAb  = \left(\int_{\Omx}\calA^2\dx\right)^{\frac{1}{2}},
  \label{eq:1265}\\[0.5em]
  &\calB(\xs)  = \left(\int_{\Omv}\left|\frac{\partial\hs}{\partial\vs}\right|^2\expAW\ddv\right)^{\frac{1}{2}},\qquad
  \calBb  = \left(\int_{\Omx}\calB^2\dx\right)^{\frac{1}{2}}.
  \label{eq:1270}
\end{align}
%%
%\begin{align}
%  \calAb 
%  = \left(\int_{\Omega} \hs^2\expAW   \ddv\dx\right)^{\frac{1}{2}}
%  = \left(\int_{\Omega} \fs^2e^{\vs^2} \ddv\dx\right)^{\frac{1}{2}}
%  = \left(\int_{\Omx} \sum_{n=0} \Big(\gamma_n \Cs_{n}(\xs,t)\Big)^2 \dx\right)^{\frac{1}{2}}
%  \label{eq:1228}
%\end{align}
%%
%\begin{align}
%  \Ls\fs     &= \frac{1}{2}\frac{\partial\fs}{\partial\vs}-\vs\fs = \frac{1}{2}\frac{\partial\hs}{\partial\vs}\expAW, \label{eq:1230}\\[0.5em]
%  \Lst\fs    &= \frac{\partial\fs}{\partial\vs},                                                                      \label{eq:1235}\\[0.5em]
%  \Lst\Ls\fs &= \frac{1}{2}\frac{\partial^2\fs}{\partial\vs^2} -\vs\frac{\partial\fs}{\partial\vs} -\fs.               \label{eq:1240}
%\end{align}
%%
Then, we rewrite problem \eqref{eq:1220} in weak form.
To this end, we multiply \eqref{eq:1220} by the test function $\phi$,
integrate over domain $\Omega=\Omx\times\Omv$, and define the bilinear form:
\begin{align}
  \Bs(\hs,\phi)
  &= \int_{\Omega}\hs\phi\expAW\ddv\dx 
  %%+ \frac{\Dt}{2}\int_{\Omega}\vs\frac{\partial}{\partial x}\phi\ddv\dx
  + \frac{\Dt}{2}\int_{\Omega}\phi\vs\left(\frac{\partial\hs}{\partial x}\right)\expAW\ddv\dx
  \nonumber\\[0.5em]
  &\quad
  - \frac{\Dt}{4}\int_{\Omx}\big( \Es^{j}+\Es^{j-1} \big)\left[\int_{\Omv}\frac{\partial\big(\hs\expAW\big)}{\partial\vs}\phi\ddv\right]\dx %% \nonumber\\[0.5em]
  %&\quad 
  + \frac{\nu \Dt}{2^{k+1}}\int_{\Omega}\frac{\partial^k\hs}{\partial\vs^k}\frac{\partial^k\phi}{\partial\vs^k}\expAW\ddv\dx,
  \label{eq:1280}
\end{align}
where the last term is obtained after successive integration by parts
as done in~\eqref{eq:0067} and using formula~\eqref{eq:0025} for
$\Ls^{(k)}\fs$. 
Now, we consider the problem of finding $\fs=\hs\expAW$ such that:
\begin{align}
  \Bs(\hs,\phi) = \int_{\Omega}\gs^{j-1}\phi\ddv\dx,
  \label{eq:1282}
\end{align}
for every $\phi$. 
Both $h$ and $\phi$ will be represented as a suitable expansion
(finite or infinite) of Hermite polynomials.  
We skip the details concerning the formulation in the proper
functional spaces, since this aspect is not relevant for the analysis
we are carrying out in this paper.

We want the bilinear form $\Bs$ to be positive definite.
First, we discuss the case $k=1$, and note that the last integral term
in~\eqref{eq:1280} can be transformed as follows
\begin{align}
  \frac{\nu\Dt}{4}\int_{\Omega}\frac{\partial\hs}{\partial\vs}\frac{\partial\hs}{\partial\vs}\expAW\ddv\dx
  %%= \frac{\nu\Dt}{4}\int_{\Omx}\left(\int_{\Omv}\frac{\partial}{\partial\vs}\left(\frac{\partial\hs}{\partial\vs}\expAW\right)\hs\ddv\right)\dx
  = \frac{\nu\Dt}{4}\calBb^2.
  \label{eq:1245}
\end{align}
%%
%%
%Recalling again that $\fs=\hs\expAW$, we use the test function
%$\phi=h$ and integrate by parts:
%\begin{align}
%  \int_{\Omv}\Lst\big(\Ls\fs\big)\hs\ddv
%  =  \frac{1}{2}\int_{\Omv}\frac{\partial}{\partial\vs}\left(\frac{\partial\hs}{\partial\vs}\expAW\right)\hs\ddv
%  = -\frac{1}{2}\calB^2.
%  \label{eq:1245}
%\end{align}
%%
%%
In this way, we get:
\begin{align}
  \Bs(\hs,\hs)
  &= \calAb^2 
  + \int_{\Omv}\frac{\Dt}{2}\vs\underbrace{\left(\frac{1}{2}\int_{\Omx}\frac{\partial\hs^2}{\partial x}\dx\right)}_{=0}\expAW\ddv
  \nonumber\\[0.5em]
  &-\frac{\Dt}{4}\int_{\Omx}\big( \Es^{j}+\Es^{j-1} \big)\left[\int_{\Omv}\frac{\partial}{\partial\vs}\left(\hs\expAW\right)\hs\ddv\right]\dx
  +\nu\frac{\Dt}{4}\calBb^2,
  \label{eq:1255}
\end{align}
where we noted that the integral of $\partial\hs^2\slash{\partial
    x}$ over $\Omx$ is zero because we assumed periodicity in
  space.
We successively integrate by parts the third term on the right:
\begin{align}
  &-\frac{\Dt}{4}\int_{\Omx}\big( \Es^{j}+\Es^{j-1} \big)\left[\int_{\Omv}\frac{\partial}{\partial\vs}\left(\hs\expAW\right)\hs\ddv\right]\dx    
%	\nonumber\\[0.5em]
 = \frac{\Dt}{4} \int_{\Omx}\big( \Es^{j}+\Es^{j-1} \big)\left[\int_{\Omv}\hs\frac{\partial\hs}{\partial\vs}\expAW\ddv\right]\dx.
% \nonumber\\[0.5em]
%   &\qquad\qquad= \frac{\Dt}{4}\int_{\Omx}\big( \Es^{j}+\Es^{j-1} \big)\left[\int_{\Omv}\hs^2v\expAW\ddv\right]\dx.
  \label{eq:1260}
\end{align}
Let us now define:
\begin{align}
  \calM = \max_{x\in\Omx}\ABS{ \Es^{j}+\Es^{j-1} }.
  \label{eq:1275}
\end{align}
Since $\calM$ depends on $\Es^{j}$ (and, consequently, on $\fs^{j}$
through~\eqref{eq:1216}, we may assume that for $\Delta t$
sufficiently small, $\Es^{j}\approx\Es^{j-1}$.
Thus, $\calM\approx 2\max_{\xs\in\Omx}\ABS{\Es^{j-1}}$.
This makes the following evaluation of $\Delta t$ practically possible
(see the estimate in~\eqref{eq:1387} below).
%% 
% This is a little trick to get out of the impasse of having to handle
% the nonlinear version.
%% 
%% We will not be able to do better in this paper, due to the extreme
%% complexity of the whole problem.

\smallskip
%%
%% Furthermore,  an 
We estimate~\eqref{eq:1260} by applying the Schwartz and Young
inequalities as follows:
\begin{align}
  &-\frac{\Dt}{4}\int_{\Omx}\big( \Es^{j}+\Es^{j-1} \big)\left[\int_{\Omv}\hs\frac{\partial\hs}{\partial\vs}\expAW\ddv\right]\dx
  \nonumber\\[0.5em]&\qquad\qquad
  \geq 
  -\frac{\Dt}{4} 
  \int_{\Omx}\abs{ \Es^{j}+\Es^{j-1} }
  \left[ \int_{\Omv}h^2 \expAW\ddv \right]^{\frac{1}{2}}
  \left[ \int_{\Omv}\left(\frac{\partial\hs}{\partial\vs}\right)^2\expAW\ddv \right]^{\frac{1}{2}}
  \,\dx
  \nonumber\\[0.5em]&\qquad\qquad
  \geq  
  -\frac{\Dt}{4} \calM \int_{\Omx}\calA\calB\dx
  %\nonumber\\[0.5em]&\qquad\qquad
  \geq 
  -\frac{\Dt}{4} \calM \int_{\Omx}\left( \frac{\sigma}{2}\calA^2 + \frac{1}{2\sigma}\calB^2 \right)\dx
  \nonumber\\[0.5em]&\qquad\qquad
  %\leq \frac{\Dt}{4}\max_{x\in\Omx}\abs{ \Es^{j}+\Es^{j-1} }\,\left( \frac{\sigma}{2}\int_{\Omx}\calA^2\dx + \frac{1}{2\sigma}\int_{\Omx}\calB^2\dx \right)
	%  \nonumber\\[0.5em]&\qquad\qquad
  =-\frac{\Dt}{4}\calM\,\left( \frac{\sigma}{2}\calAb^2 + \frac{1}{2\sigma}\calBb^2 \right),
\label{eq:1295}
\end{align}
where $\sigma >0$ is an arbitrary parameter.
Using this estimate in~\eqref{eq:1255}, we find the inequality
\begin{align}
  \Bs(\hs,\hs)
  \geq \calAb^2 +\nu\frac{\Dt}{4}\calBb^2
  -\frac{\Dt}{4}\calM\,\left( \frac{\sigma}{2}\calAb^2 + \frac{1}{2\sigma}\calBb^2 \right).
   \label{eq:1295a}
\end{align}
To derive sufficient conditions for the positivity of the bilinear form
$\Bs(\cdot,\cdot)$, i.e., $\Bs(\hs,\hs)\geq0$, we can proceed in
different ways.
First, for every strictly positive quantity $\sigma$, we can impose that 
\begin{align}
  \frac{\Dt}{4}\calM\,\left( \frac{\sigma}{2}\calAb^2 + \frac{1}{2\sigma}\calBb^2 \right)
  \leq \calAb^2 + \nu\frac{\Dt}{4}\calBb^2.
  \label{eq:1305}
\end{align}
%or, after neglecting the surely positive term with the numerical diffusion parameter $\nu$, that
%\begin{align}
%  \calAb^2 -\frac{\Dt}{4}\calM\,\left( \frac{\sigma}{2}\calAb^2 + \frac{1}{2\sigma}\calBb^2 \right)\geq 0.
%  \label{eq:1305:b}
%\end{align}
%%%
%In the first case, we expect to find a sufficient condition that relates 
%$\nu$, $\N$, $\Delta\ts$, and $\calM$.
%%
%In the second case, we expect to find a sufficient condition that relates only 
%$\N$, $\Delta\ts$, and $\calM$ since $\nu$ has been removed from the inequality.
%
To this end, we note that:
%To elaborate on~\eqref{eq:1305}, we note that:
\begin{align*}
    \frac{\sigma}{2}\calAb^2 + \frac{1}{2\sigma}\calBb^2
    =\frac{\sigma}{2}
    \left( \calAb^2 + \frac{1}{\sigma^2}\calBb^2 \right).
\end{align*}
Comparing the expression above with the right-hand side of inequality~\eqref{eq:1305}, 
suggests us to set $1\slash{\sigma^2}=\nu\Dt\slash{4}$, or, equivalently that
$\sigma = 2\slash{\sqrt{\nu\Delta t}}$.
We set this value of $\sigma$ back into inequality~\eqref{eq:1305} to find that
\begin{align}
  \frac{\Dt}{4}\calM\,
  \frac{\sigma}{2}\calAb^2 + \frac{1}{2\sigma}\calBb^2
  = \frac{\Dt\calM}{4\sqrt{\nu\Delta t}}\left(\calAb^2 + \nu\frac{\Delta t}{4}\calBb^2\right)
  \leq \left(\calAb^2 + \nu\frac{\Delta t}{4}\calBb^2\right),
  \label{eq:1385a}
\end{align}
from which we immediately have the condition:
\begin{align}
  \frac{\Dt\calM}{4\sqrt{\nu\Delta t}}\leq 1,
\label{eq:1385bb}
\end{align}
that we can rewrite as
\begin{align}
  \Delta t\leq \frac{16\nu}{\calM^2},
  \label{eq:1387}
\end{align}
after renormalizing the factor $\sqrt{\Dt}$ in the denominator of~\eqref{eq:1385bb} 
and squaring the resulting inequality.
Such a constraint on $\Delta t$ constitutes a sufficient condition to
realize the invertibility of problem~\eqref{eq:1220} for
$k=1$. 
Unfortunately, we are unable to provide a similar result in the case
when $k>1$.  
The problem is that inequality~\eqref{eq:1305} becomes of the form:
\begin{align}
  \frac{\Dt}{4}\calM\,\left( \frac{\sigma}{2}\calAb^2 + \frac{1}{2\sigma}\calBb^2 \right)
  \leq \calAb^2 + \frac{\nu \Dt}{2^{k+1}}\int_{\Omega}\left(\frac{\partial^k\hs}{\partial\vs^k}\right)^2\expAW\ddv\dx .
  \label{eq:1312}
\end{align}
We can bound $\calBb$, that contains only first derivatives, by an
expression containing higher order derivatives, only if a certain
number of low modes of $h$ is set to zero.
This is certainly not consistent with the freedom we would like to
leave to these coefficients.

\smallskip
To recover an alternative estimate of the time step $\Delta\ts$ that does not
involve the diffusion parameter $\nu$,
we suppose that $h$ is a linear combination
of a finite number of Hermite polynomials.
In practice, $h$ is going to be a polynomial of degree less than or equal
to $N$. 
In this situation, we can rely on the inverse type inequality:
\begin{align}
  \calBb \ \leq \ \sqrt{2N} \ \calAb,
\end{align} 
which is easily deducible from~\eqref{eq:2502a}.
Thus, to control the last term at the end of~\eqref{eq:1295a} 
we proceed by writing:
\begin{align}
  -\frac{\Dt}{4}\calM\,\left( \frac{\sigma}{2}\calAb^2 + \frac{1}{2\sigma}\calBb^2 \right)
  \geq -\frac{\Dt}{4}\calM\,\left( \frac{\sigma}{2} + \frac{N}{\sigma}\right)\calAb^2
  =    -\frac{\Dt}{4}\calM \sqrt{2N}\,\calAb^2,
  \label{eq:1296}
\end{align}
where we noticed that the absolute value of the term in the middle is 
minimized by the choice $\sigma =\sqrt{2N}$.
%% and $1/2\sqrt{N\slash{2}}+N\sqrt{2\slash{N}}\leq 2\sqrt{2N}$
%%
In this way, the positivity of the bilinear form is realized by
requiring that the last term in~\eqref{eq:1296} is less than
$\calAb^2+\frac14 \nu \Delta t \calBb^2$.
This is true by choosing:
\begin{align}
  \Delta t\leq \frac{4}{\calM \sqrt{2N}},
  \label{eq:1388}
\end{align}
and, now, the bound on $\Delta\ts$ depends on $N$ but not on $\nu$.
Moreover, this calculation does not involve any explicit
expression from the Lenard-Bernestein diffusion operators 
on the right-hand side of~\eqref{eq:1280} since this term
was just eliminated because of its positivity for $\phi =h$.
This means that this time the relation between $N$, $\Delta\ts$, and $\calM$
holds for any value of $k\geq 1$.

\smallskip
We can make further considerations by putting together
inequalities~\eqref{eq:1387} and ~\eqref{eq:1388}.  
If $\Delta t$ is chosen in order to be consistent with both of them, we
get:
\begin{align}
  \Delta t \approx \frac{16\nu}{\calM^2}
  \quad\textrm{and}\quad
  \Delta t\approx \frac{4}{\calM \sqrt{2N}}
  \quad\textrm{imply~that}\quad
  \ \nu\approx\frac{\calM}{4\sqrt{2N}}. %% \quad N\nu\Delta t\approx \frac12
  \label{eq:1587}
\end{align}
Similarly, setting $1\slash{\calM}=\sqrt{2N}\Delta t\slash{4}$ in
$\Delta t\approx 16\nu\slash{\calM^2}$ above implies that
\begin{align}
  \Delta t \approx 16\nu \frac{1}{\calM^2} = 16\nu\frac{2N\Delta t^2}{16} = 2N\nu\Delta t^2,
\end{align}
from which we derive the relation
\begin{align}
  \nu N\Delta t\approx \frac12.
\end{align}
The last relation agrees with the suggestion, made in the previous section, that the product $N\nu\Delta t$
should be of order of the unity.

%At this point, we observe that $\vert v\vert \leq 1+v^2$. This
%allows us to use the inequality~\eqref{eq:2600}. Recalling~\eqref{eq:1275}, we bound
%the last termi in~\eqref{eq:1260} as follows: 
%\begin{align}
%\left\vert\frac{\Dt}{4}\int_{\Omx}\big( \Es^{j}+\Es^{j-1} \big)\left[\int_{\Omv}\hs^2v\expAW\ddv\right]\dx
%\right\vert \leq \frac{\Dt}{4}\calM \int_{\Omega}\hs^2(1+v^2)\expAW\ddv\dx
%  \label{eq:1266}
%\end{align}

\smallskip
We can say something more if the electric field is treated explicitly, i.e.: $(E^j+E^{j-1})/2\approx E^{j-1}$,
with $\calM  =2\max_{x\in{\Omx}}\vert E^{j-1}\vert$. 
The maximum norm can be bounded through the first derivative.
This is done in the following way:
\begin{align}
  \calM^2 \leq \vert\Omx\vert \int_{\Omx} \left(\frac{\partial E^{j-1}}{\partial x}\right)^2 dx =
  \vert \Omx\vert \int_{\Omx} \left( 2- \int_{\Omv} f^{j-1} dv\right)^2 dx,
  \label{eq:1485}
\end{align}
where $\vert\Omx\vert$ denotes the measure of $\Omx$. 
Next, we use a standard inequality and the Schwartz inequality to obtain:
\begin{align}
  \frac{\calM^2}{ \vert\Omx\vert }
  &\leq \int_{\Omx}\left( 2-\int_{\Omv}\fs^{j-1} dv\right)^2 dx
  \leq 2\int_{\Omx}\left[ 4+ \left(\int_{\Omv} f^{j-1} dv\right)^2 \right]dx
  \nonumber\\[0.5em]
  &\leq 8\vert\Omx\vert + \int_{\Omx}
  \left[ 
    \int_{\Omv} (f^{j-1})^2 e^{ v^2} dv 
    \int_{\Omv}             e^{-v^2} dv
  \right]dx
  \nonumber\\[0.5em]
  &\leq 8 \vert \Omx\vert + \sqrt{\pi} \int_{\Omx} \int_{\Omv} (h^{j-1})^2 e^{-v^2} dv dx
  =     8 \vert \Omx\vert + \sqrt{\pi} \cal{H},
  \label{eq:1414}
\end{align}
where we denoted the last integral by $\cal{H}$.
Then, we consider again inequality~\eqref{eq:1295a} from which we remove the
nonnegative term $\nu\Dt\calBb^2\slash{4}$ to obtain a sufficient condition 
that is independent of $\nu$.
Using~\eqref{eq:1414} in the right-hand side of~\eqref{eq:1296}, we end up with
\begin{align}
  \frac{\Delta t}{4}\vert\Omx\vert^{1/2}\left( 8\vert\Omx\vert + \sqrt{\pi}\cal{H} \right)^{1/2}
  \sqrt{2N}\calAb^2 \leq \calAb^2,
  \label{eq:1385b}
\end{align}
which implies
\begin{align}
\Delta t \vert\Omx\vert^{1/2} \left( 8\vert\Omx\vert + \sqrt{\pi} \cal{H} \right)^{1/2}
  \leq \frac{4}{\sqrt{2N}},
  \label{eq:1385c}
\end{align}
This last condition is substantially similar to~\eqref{eq:1388}. 
However, this derivation implies that having a knowledge of either
$\calM$ or $\cal{H}$ at the step $j-1$, we have an idea on how to set
up the new time-step for the successive iteration.

\smallskip
In the final part of our study, we put together what we have learned in
the previous sections, and investigate the interplay between
time stability and conservation properties. 
We consider, first, the conservation of the mass, which is the zero-th 
order moment of the Vlasov distribution function $\fs$.
After discretization in time, we assume that $f^j$ is expanded on the
Hermite functions' basis:
\begin{align}
  f^j(x,v) = \sum_{n=0}^\infty C_n^{\star,j}(x) \psi_n(v).
  \label{eq:1488}
\end{align}
The variational formulation for the expansion coefficients $C_n^{\star,j}$ is obtained by substituting~\eqref{eq:1488}
in~\eqref{eq:1215}, multiplying by the test function $\psi^m$ and integrating on $\Omega=\Omx\times\Omv$:
%
%We simplify the integrals involving $\psi_n\psi^m$ and $\partial\psi_n\slash{\partial\vs}\,\psi^m$
%by using the orthogonality properties of the Hermite polynomials, and we find:
\begin{align}
  &\sum_{n=0}^\infty\left[\int_\Omega \frac{C_n^{\star,j}-C_n^{\star,j-1}}{\Delta t}\psi_n\psi^m dxdv \right]+
  \sum_{n=0}^\infty\left[\int_\Omega \frac{\partial}{\partial x}\left(\frac{C_n^{\star,j}+C_n^{\star,j-1}}{2}
    \right)v\psi_n\psi^m dxdv\right]  
  \nonumber\\[0.5em]
  &\qquad
  -\sum_{n=0}^\infty\left[\int_{\Omega}\frac{\Es^{j}+\Es^{j-1}}{2}\ \frac{C_n^{\star,j}+C_n^{\star,j-1}}{2}
    \frac{\partial \psi_n}{\partial v}\psi^m dxdv \right] 
  \nonumber\\[0.5em]
  &\qquad
  +(-1)^k\nu\sum_{n=0}^\infty\left[\int_\Omega \Lst^{(k)}\Ls^{(k)}
    \left(\frac{C_n^{\star,j}+C_n^{\star,j-1}}{2}\right)\psi_n\psi^m dxdv \right]=0.
  \label{eq:1460}
\end{align}
We separate the integration with respect to $x$ from that
with respect to $v$, obtaining:
\begin{align}
  &\int_{\Omx} \frac{C_m^{\star,j}-C_m^{\star,j-1}}{\Delta t} dx +
  \sum_{n=0}^\infty\left[\int_{\Omx} \frac{\partial}{\partial x}\left(\frac{C_n^{\star,j}+C_n^{\star,j-1}}{2}\right)dx   
    \int_{\Omv} v\psi_n\psi^m dv \right]
  \nonumber\\[0.5em]
  &\qquad
  + \frac{\gamma_m}{\gamma_{m+1}}\int_{\Omx}\frac{\Es^{j}+\Es^{j-1}}{2}\ \frac{C_{m-1}^{\star,j}+C_{m-1}^{\star,j-1}}{2}dx  
  \nonumber\\[0.5em]
  &\qquad
  -m(m-1)\cdots (m-k+1)\nu\int_{\Omx} 
  \frac{C_m^{\star,j}+C_m^{\star,j-1}}{2}\dx =0.
  \label{eq:1461}
\end{align}
We further note that, due to the periodic boundary conditions, the
integral in the variable $x$ of the second term is zero.
In terms of the coefficients in the Hermite polynomial
basis, Eq. (\ref{eq:1461}) becomes:
\begin{align}
  &\int_{\Omx} \frac{C_n^{j}-C_n^{j-1}}{\Delta t} dx 
  +\sqrt{\frac{n+1}{n}}\int_{\Omx}\frac{\Es^{j}+\Es^{j-1}}{2}\ \frac{C_{n-1}^{j}+C_{n-1}^{j-1}}{2}dx  
  \nonumber\\[0.5em]
  &\qquad
  -n(n-1)\cdots (n-k+1)\nu\int_{\Omx} 
  \frac{C_n^{j}+C_n^{j-1}}{2}\dx =0.
  \label{eq:1461a}
\end{align}
This system of equations is coupled with~\eqref{eq:1216}.
As a consequence of the orthogonality, we have:
\begin{align}
    \int_{\Omv}\fs^{j}\ddv = \sum_{n=0}^\infty \int _{\Omv}C_n^{j}H_n \expAW = \sqrt{\pi}\ C^j_0.
  \label{eq:1226}
\end{align}
Thus, the discretized Poisson equation takes the form:
\begin{align}
   \frac{\partial \Es^{j}}{\partial x} = 1 - \sqrt{\pi}\ C^j_0.
  \label{eq:1227}
\end{align}
By integrating this last relation with respect to $x$ and using the
boundary conditions for $E^j$, we discover that $\int_{\Omx}C^j_0 dx$
is constant for all $j\geq 0$.
This condition is maintained by the scheme~\eqref{eq:1461a}, whatever
is $k\geq 1$.
More in general, conservation of momenta $\int_\Omega v^m f^j dxdv$,
$j\geq 0$, is guaranteed up to $m\leq k-1$. 
This corresponds to the generalization for arbitrary k of the
conservation properties that were proven in
Ref. \cite{Camporeale-Delzanno-Bergen-Moulton:2016} for k=3.
\section{Conclusion}
\label{sec:conclusion}

We investigated the role of Lenard-Bernstein-like pseudo-collisional
operators in conjunction with spectral approximations of the Vlasov
equation for a collisionless plasma in the electrostatic limit.
In particular, we analyzed the spectral approximation of some
one-dimensional, simplified model problems based on different families
of Hermite basis functions using the symmetric and the asymmetric
formulations.
In the asymmetric case, we were able to prove the absolute stability
in time in an $\LTWO$-weighted norm, a problem that has been
unresolved for many years.
The results have partially been extended to the case of the full
Vlasov-Poisson model.

%% ----------------
%% Acknowledgements
%% ----------------

\section*{Acknowledgements}
This work was supported by the LDRD program of Los
Alamos National Laboratory under project number 20170207ER.
Los Alamos National Laboratory is operated by Triad National Security,
LLC, for the National Nuclear Security Administration of
U.S. Department of Energy (Contract No. 89233218CNA000001).
The authors are affiliated to the Italian
Istituto Nazionale di Alta Matematica (INdAM).
This manuscript has no associated data.

%The first and second author (DF and GM) are affiliated to the Italian
%Istituto Nazionale di Alta Matematica (INdAM).

%% Bibliography
%% \bibliographystyle{plain}
%% \bibliography{ref}

\end{document}